\documentclass[11pt]{article}
\usepackage{mathrsfs}
\usepackage{amsmath,amsfonts,amssymb,rotating,amsthm}
\usepackage{hyperref}
\usepackage{array}
\usepackage{url,booktabs}
\usepackage{amscd, amsfonts, epsfig, graphicx, colordvi,verbatim,ifpdf, pifont,graphics,cite}
\usepackage{psfrag,eepic,color}
\newtheorem{theo}{Theorem}[section]

\newtheorem{conj}[theo]{Conjecture}
\theoremstyle{remark}

\numberwithin{equation}{section}

\newdimen\Squaresize \Squaresize=11pt
\newdimen\Thickness \Thickness=0.7pt
\def\Square#1{\hbox{\vrule width \Thickness
   \vbox to \Squaresize{\hrule height \Thickness\vss
    \hbox to \Squaresize{\hss#1\hss}
   \vss\hrule height\Thickness}
\unskip\vrule width \Thickness} \kern-\Thickness}

\def\Vsquare#1{\vbox{\Square{$#1$}}\kern-\Thickness}

\def\moins{\raise 1pt\hbox{{$\scriptstyle -$}}}
\newcommand{\ospt}{\mathop{\mathrm{ospt}}\nolimits}

\newcommand{\spt}{\mathop{\mathrm{spt}}\nolimits}
\newcommand{\Spt}{\mathop{\mathrm{Spt}}\nolimits}
\newcommand{\M}{\mathop{\mathrm{M2spt}}\nolimits}
\newcommand{\sptcrank}{\mathop{\mathrm{sptcrank}}\nolimits}

\newtheorem{thm}{Theorem}[section]

\newtheorem{conje}[thm]{Conjecture}
\newtheorem{defi}[thm]{Definition}

\newtheorem{core}[thm]{Corollary}

\allowdisplaybreaks

\begin{document}
\begin{center}
{\Large \bf  The spt-Function of Andrews }
\end{center}

\begin{center}
{\large William Y.C. Chen}\\[9pt]

Center for Applied Mathematics\\
Tianjin University\\
Tianjin 300072, P. R. China\\[6pt]

and

Center for Combinatorics, LPMC\\
Nankai University\\
 Tianjin 300071, P. R. China\\[6pt]
Email: {\tt
chen@nankai.edu.cn\tt }
\end{center}

\vspace{0.3cm} \noindent{\bf Abstract.}
The spt-function $\spt(n)$ was introduced by Andrews as
the weighted counting of partitions of $n$ with respect to
the number of occurrences of the smallest part. Andrews showed that $\spt(5n+4)\equiv 0 \pmod{5}, \spt(7n+5)\equiv 0 \pmod{7}$ and $\spt(13n+6)\equiv 0 \pmod{13}$.
Since then, congruences of $\spt(n)$ have been extensively studied. Folsom and Ono
obtained congruences   of $\spt(n)$ mod  $2$ and $3$. They also showed that
the generating function of $\spt(n)$ mod 3 is related to a weight 3/2 Hecke eigenform with Nebentypus.
Combinatorial interpretations of  congruences of $\spt(n)$ mod $5$ and $7$ have been found  by
Andrews, Garvan and Liang by introducing the spt-crank of a vector partition.
Chen, Ji and Zang showed that the
 set of partitions counted by $\spt(5n+4)$
(or $\spt(7n+5)$) can be divided into five (or seven) equinumerous classes
 according to the spt-crank of a doubly marked partition.
 Let $N_S(m,n)$  denote the net number of $S$-partitions of $n$ with spt-crank $m$.
Andrews, Dyson and Rhoades conjectured that $\{N_S(m,n)\}_m$ is unimodal for any $n$.
 Chen, Ji and Zang gave a constructive proof of this conjecture.
In this survey, we summarize developments on congruence properties of $\spt(n)$
established by Andrews, Bringmann, Folsom, Garvan, Lovejoy and Ono et al.,
as well as  their  combinatorial interpretations. Generalizations and variations of the spt-function are also  discussed.
We further give an overview of asymptotic formulas of $\spt(n)$ obtained by Ahlgren, Andersen and Rhoades et al.
We conclude with some conjectures  on inequalities on  $\spt(n)$, which are reminiscent of inequalities on $p(n)$  due to DeSalvo and Pak, and Bessenrodt and Ono. Furthermore, we observe that, beyond the log-concavity,
$p(n)$ and $\spt(n)$ satisfy higher order inequalities based on polynomials arising in the invariant theory of
binary forms. In particular, we conjecture that the higher order Tur\'{a}n inequality $4(a_n^2-a_{n-1}a_{n+1})(a_{n+1}^2-a_{n}a_{n+2})-(a_na_{n+1}-a_{n-1}a_{n+2})^2>0$  holds for $p(n)$ when $n\geq 95$ and for $\spt(n)$ when $n\geq 108$.

\section{Introduction}

Andrews \cite{Andrews-2008} introduced the spt-function $\spt(n)$ as the
weighted counting of partitions with respect to the
number of occurrences of the smallest part and he discovered
that the spt-function bears striking resemblance to the classical
partition function $p(n)$. Since then, the spt-function has drawn
 much attention and has been extensively studied. In this survey, we shall summarize
 developments on the spt-function including congruence properties derived from
  $q$-identities and modular forms, along with their  combinatorial interpretations, as well as generalizations, variations and asymptotic properties.
 For the background on partitions, we refer to \cite{Andrews-1979, Andrews-1976, Andrews-Eriksson-2004},
 and for the background on modular forms, we refer to
  \cite{Apostol-1990, Diamond-Shurman-2005, Koblitz-1993, Ono-2004}.

The spt-function $\spt(n)$, called the smallest part function, is defined to be the total number of smallest parts in all partitions of $n$. More precisely, for a partition $\lambda$ of $n$, we use $n_s(\lambda)$ to denote
the number of occurrences of the smallest part in $\lambda$. Let $P(n)$ denote the set of partitions of $n$, then
\begin{equation}\label{defi-spt-1}
\spt(n)=\sum_{\lambda \in P(n)}n_s(\lambda).
\end{equation}
For example, for $n=4$, we have $\spt(4)=10$. Partitions in $
P(4)$ and the values of $n_s(\lambda)$ are
listed below:
\[\begin{array}{c|c|c|c|c|c}
\lambda \in  P(4)& (4)&(3,1)&(2,2)&(2,1,1)&(1,1,1,1)\\[5pt] \hline
  n_s(\lambda) &1&1&2&2&4
\end{array}
\]
 The spt-function $\spt(n)$ can also be interpreted by marked partitions,
   see Andrews, Dyson and Rhoades \cite{Andrews-Dyson-Rhoades-2013}. A marked partition of $n$  is meant to be  a pair $(\lambda,k)$, where $\lambda=(\lambda_1,\lambda_2,\ldots,\lambda_l)$ is an ordinary  partition of $n$ and $k$ is an integer identifying one of its smallest parts. If $\lambda_k$ is the identified  smallest part of $\lambda$, we then use $(\lambda,k)$
   to denote this marked partition.   For example, there are ten marked partitions of $4$.
  \[\small \begin{array}{lllll}
((4),1),&  ((3,1), 2),& ((2,2),1),&((2,2),2),&((2,1,1),2),\\[5pt]
((2,1,1),3),&((1,1,1,1),1),&((1,1,1,1),2),&((1,1,1,1),3),&((1,1,1,1),4).
\end{array}
\]

 Using the definition \eqref{defi-spt-1}, it is easy to derive  the following
 generating function,  see Andrews \cite{Andrews-2008},
\begin{equation}\label{gf-spt}
\sum_{n=1}^\infty\spt(n)q^n  =\sum_{n=1}^\infty \frac{q^n}{(1-q^n)^2(q^{n+1};q)_\infty}.
\end{equation}

Here we have adopted the common   notation \cite{Andrews-1976}:
\begin{align*}
(a;q)_\infty=\prod_{n=0}^\infty(1-aq^n) \quad \text{and} \quad
(a;q)_n&=\frac{(a;q)_\infty}{(aq^n;q)_\infty}.
\end{align*}

The spt-function is closely related to the rank and the crank of a partition. Recall that the rank of a partition was
    introduced by Dyson \cite{Dyson-1944} as the largest part of the partition minus the number of parts.  The crank of a partition was defined by Andrews and Garvan \cite{Andrews-Garvan-1988} as the largest part if the partition contains no ones,  otherwise as the number of parts larger than the number of ones minus the number of ones.   For $n\geq 1$, let $N(m,n)$ denote the number of partitions of $n$ with rank $m$,
and for $n>1$, let $M(m,n)$ denote the number of partitions of $n$ with crank $m$. For $n=1$, set
\[M(0,1)=-1,\,M(1,1)=M(-1,1)=1,\,\]
and for $n=1$ and  $m\neq -1,0,1$, set
\[M(m,1)=0.\]

 Atkin and Garvan \cite{Atkin-Garvan-2011} defined the $k$-th moment $N_k(n)$  of ranks   as
  \begin{eqnarray}\label{equ-def-nkmn}
N_k(n)&=& \sum_{m=-\infty}^{\infty}m^kN(m,n),
\end{eqnarray}
and the $k$-th  moment   $M_k(n)$   of cranks as
\begin{equation*}
M_k(n)=\sum_{m=-\infty}^\infty m^k M(m,n)\label{equ-def-mk}.
\end{equation*}
It is worth mentioning that Atkin and Garvan  \cite{Atkin-Garvan-2011}  showed that the generating functions of the moments of  cranks  are related to quasimodular forms. Bringmann, Garvan and Mahlburg \cite{Bringmann-Garvan-Mahlburg-2009} showed that  the generating functions of the moments of ranks  are related to quasimock theta functions. Asymptotic formulas for the moments of  ranks and cranks  were derived by Bringmann, Mahlburg and Rhoades \cite{ Bringmann-Mahlburg-Rhoades-2014}.

 Based on the generating function \eqref{gf-spt} and Watson's $q$-analog of Whipple's theorem \cite[p.\,43, eq.\,(2.5.1)]{Gasper-Rahman-2004},
 Andrews \cite{Andrews-2008} showed that the spt-function can be expressed in terms of the second moment $N_2(n)$  of ranks introduced by  Atkin and Garvan \cite{Atkin-Garvan-2011},
\begin{equation}\label{spt-moments}
\spt(n)=np(n)-\frac{1}{2}N_2(n).
\end{equation}
Ji  \cite{Ji-2008} found a combinatorial proof of \eqref{spt-moments}  using rooted partitions.

By means of a   relation due to Dyson \cite{Dyson-1989}, namely,
\begin{equation}\label{m2-2npn-dyson}
M_2(n)=2np(n),
\end{equation}
Garvan \cite{Garvan-2010} observed that the  expression
\begin{equation}\label{spt-def-eq}
\spt(n)=\frac{1}{2}M_2(n)-\frac{1}{2}N_2(n)
\end{equation}
 implies that $M_2(n)>N_2(n)$ for $n\geq 1$. In general, he
 conjectured and later proved that $M_{2k}(n)>N_{2k}(n)$ for $k\geq 1$ and $n\geq 1$, see \cite{Garvan-2010, Garvan-2011}.

In view of  the relation \eqref{spt-moments} and identities
on  refinements of $N(m,n)$ established by Atkin and Swinnerton-Dyer \cite{Atkin-Swinnerton-Dyer-1954} and O'Brien  \cite{O'Brien-1965}, Andrews proved that
 $\spt(n)$ satisfies  congruences mod $5$, $7$ and $13$ which are
 reminiscent of Ramanujan's congruences for $p(n)$. Let  $\ell$ be a prime. A Ramanujan congruence modulo $\ell$ for the sequence  $\{a(n)\}_{n\geq 0}$ means a congruence of the form
 \[a(\ell n+\beta)\equiv 0 \pmod{\ell}\]
 for  all nonnegative integers $n$ and a fixed integer $\beta$.

  Ramanujan \cite{Ramanujan-1919} discovered   the following   congruences for $p(n)$,
 \begin{eqnarray}
p(5n+4)&\equiv& 0 \pmod 5, \label{pcon-5-16}\\[3pt]
p(7n+5) &\equiv& 0 \pmod 7,\label{pcon-7-17}\\[3pt]
p(11n+6) &\equiv& 0 \pmod{11},\label{pcon-7-11}
\end{eqnarray}
 and proclaimed that ``it appears that there are no equally simple properties for any moduli involving primes other than these three (i.e. $\ell=5,7,11$).'' See also Berndt \cite[p.\,27]{Berndt-2006}.

 Elementary proofs of the congruences \eqref{pcon-5-16} and \eqref{pcon-7-17} were given by Ramanujan \cite{Ramanujan-1919} and  an elementary proof of the  congruence \eqref{pcon-7-11} was given by Winquist \cite{Winquist-1969}. Alternative  proofs of   \eqref{pcon-7-11}  were found by  Berndt, Chan, Liu and Yesilyurt \cite{Berndt-Chan-Liu-Yesilyurt-2004} and Hirschhorn \cite{Hirschhorn-2013}. Recently, Paule and Radu \cite{Paule-Radu} found a recurrence relation of the generating function of $p(11n+6)$, from which
  \eqref{pcon-7-11} is an immediate consequence.   Berndt \cite{Berndt-2007}  provided simple proofs of   \eqref{pcon-5-16}--\eqref{pcon-7-11} by using Ramanujan's differential equations for the Eisenstein series.  Uniform proofs of   \eqref{pcon-5-16}--\eqref{pcon-7-11}  were found by Hirschhorn \cite{Hirschhorn-1994}.

    Concerning Ramanujan's conjecture, Kiming and Olsson \cite{Kiming-Olsson-1992}
     showed that if there exists a Ramanujan's congruence   $p(\ell n+\beta)\equiv 0\pmod{\ell}$, then $24\beta\equiv 1\pmod{\ell}$. According to this condition,  Ahlgren and Boylan \cite{Ahlgren-Boylan-2005} confirmed Ramanujan's conjecture. More precisely, they  showed that for a prime $\ell$, if there is a Ramanujan's congruence  modulo $\ell$ for $p(n)$, then it must be one of the congruences \eqref{pcon-5-16}, \eqref{pcon-7-17} and \eqref{pcon-7-11}.

  Combinatorial studies of Ramanujan's congruences of $p(n)$ go back to Dyson \cite{Dyson-1944}. He    conjectured that the rank of a partition  can be used to divide the set of partitions of $5n+4$ (or $7n+5$) into five (or seven) equinumerous classes. More precisely,
let $N(i,t,n)$ denote the number of partitions of $n$ with rank congruent to $i$ modulo $t$. Dyson \cite{Dyson-1944} conjectured that
    \begin{eqnarray}\label{Dyson-1}
    N(i,5,5n+4)&=&\frac{p(5n+4)}{5}\quad\text{for } \quad 0\leq i\leq 4,\\[3pt]
    N(i,7,7n+5)&=&\frac{p(7n+5)}{7}\quad\text{for } \quad 0\leq i\leq 6. \label{Dyson-2}
    \end{eqnarray}
    These relations were proved by Atkin and Swinnerton-Dyer \cite{Atkin-Swinnerton-Dyer-1954},
     which imply  \eqref{pcon-5-16}    and \eqref{pcon-7-17}. Dyson also pointed out that the rank of a partition cannot be used to interpret   \eqref{pcon-7-11}.    To give a combinatorial explanation  of this congruence modulo 11, Garvan \cite{Garvan-1988}  introduced the crank of a vector partition and showed that this statistic leads to interpretations of the above congruences of $p(n)$ mod 5, 7 and 11.  Andrews and Garvan \cite{Andrews-Garvan-1988} found  an equivalent definition of the crank in terms of  an ordinary  partition. For the history of the rank and the crank, see, for example, Andrews and Berndt \cite{Andrews-Berndt-2012} and Andrews and Ono \cite{Andrews-Ono-2005}.

Although Dyson's rank fails to explain Ramanujan's congruence \eqref{pcon-7-11} combinatorially, the generating functions for the rank differences  have been extensively studied. For example, the generating functions for  the rank differences $N(s, \ell, \ell n + d)- N(t, \ell, \ell n + d)$ for $\ell=2,9,11,12,13$  have been determined by Atkin and Hussain \cite{Atkin-Hussain-1958},  O'Brien   \cite{O'Brien-1965}, Lewis \cite{Lewis-1991, Lewis-1997} and Santa-Gadea \cite{Gadea-1992}.

By the relations \eqref{spt-moments}, \eqref{Dyson-1} and \eqref{Dyson-2}, Andrews \cite{Andrews-2008} showed that
    \begin{eqnarray}
 \spt(5n+4)&\equiv& 0 \pmod 5, \label{con-5}\\[3pt]
\spt(7n+5) &\equiv& 0 \pmod 7 \label{con-7}.
\end{eqnarray}
He also obtained that
\begin{equation}\label{con-13}
\spt(13n+6) \equiv 0 \pmod {13},
\end{equation}
by considering the properties of  $N(i,13,13n+6)$  due to O'Brien   \cite{O'Brien-1965}.
Let
\[r_{a,b}(d)=\sum_{n=0}^\infty(N(a,13,13n+d)-N(b,13,13n+d))q^{13n},\]
and for $1\leq i \leq 5$, and let
\[
S_i(d)=r_{(i-1),i}(d)-(7-i)r_{5,6}(d).
\]
O'Brien  \cite{O'Brien-1965} deduced that
\begin{equation}\label{identity-mod13-1}
S_1(6)+2S_2(6)-5S_5(6)\equiv 0\pmod{13}
\end{equation}
and
\begin{equation}\label{identity-mod13-2}
S_2(6)+5S_3(6)+3S_4(6)+3S_5(6)\equiv 0\pmod{13}.
\end{equation}
Employing \eqref{spt-moments}, Andrews derived an expression
for $\spt(13n+6)$ in terms of   $N(i,13, 13n+6)$ modulo $13$. Then the congruence \eqref{con-13} follows from  \eqref{identity-mod13-1} and \eqref{identity-mod13-2}.

This paper is organized as follows. In  Section 2, we recall the spt-crank of an $S$-partition defined by Andrews, Garvan and Liang, which leads to combinatorial interpretations of the congruences of the spt-function mod 5 and 7.
Motivated by a problem of Andrews, Garvan and Liang  on constructive proofs of
the congruences mod 5 and 7, Chen, Ji and Zang introduced the
  notion of a doubly marked partition and its spt-crank. Such an spt-crank
   can be used to  divide the set counted by $\spt(5n+4)$ (resp. $\spt(7n+5)$) into five (resp. seven) equinumerous classes.  The unimodality of the spt-crank and related topics are also discussed.  In Section 3, we  begin  with  Ramanujan-type congruences of $\spt(n)$ mod 11,17, 19, 29, 31 and 37 obtained by Garvan. We then consider Ramanujan-type congruences of $\spt(n)$ modulo any prime $\ell\geq 5$ due to Ono and the $\ell$-adic generalization   due to Ahlgren,  Bringmann and  Lovejoy. The congruences of $\spt(n)$ mod powers of $5$, $7$ and $13$ established by  Garvan will also be discussed. We finish this section with
   congruences of $\spt(n)$ mod $2,3$ and  powers of $2$ due to Folsom and Ono, and Garvan and Jennings-Shaffer.
    Section 4 is devoted to  generalizations and variations of the spt-function. We first recall the higher order spt-function defined by Garvan,  as a generalization of the spt-function.   We then concentrate on two generalizations of the spt-function based on the $j$-rank, given by Dixit and Yee.
     The first variation of the spt-function was defined by Andrews, Chan and Kim as the difference between  the first rank and crank moments.
   At the end of this section, we present three variations of the spt-function,
   which are restrictions of the spt-function to  three classes of  partitions.  The generating functions, combinatorial interpretations and congruences of these generalizations and variations of the spt-function will also be discussed.   In Section 5, we   summarize  asymptotic formulas of  the spt-function and its variations.   Section 6 contains some  conjectures on inequalities on $\spt(n)$,  which are analogous to  those on $p(n)$, due to DeSalvo and Pak, and Bessenrodt and Ono.  Beyond the log-concavity,
 we conjecture that $p(n)$ and $\spt(n)$ satisfy higher order inequalities induced from invariants of
binary forms. In particular, we conjecture that the higher order  Tur\'{a}n inequality
holds for both $p(n)$ and $\spt(n)$  when  $n$ is large enough.

\section{The spt-crank}

 To give combinatorial interpretations of congruences on $\spt(n)$,
 Andrews, Garvan and Liang \cite{Andrews-Garvan-Liang-2011}  introduced the spt-crank of an $S$-partition, which is analogous to Garvan's crank  of a vector partition \cite{Garvan-1988}.  They showed that the spt-crank of an $S$-partition can be used   to divide the set of $S$-partitions with signs counted by $\spt(5n+4)$ (or $\spt(7n+5)$) into five (or seven) equinumerous classes which leads to the congruences  \eqref{con-5} and \eqref{con-7}.

Andrews, Dyson and Rhoades \cite{Andrews-Dyson-Rhoades-2013} proposed
the problem of finding an equivalent definition of the spt-crank for a marked partition.   Chen, Ji and Zang \cite{Chen-Ji-Zang-2016} introduced the structure of a doubly  marked partition and established a bijection between marked partitions and  doubly marked
partitions. Then they defined the spt-crank of a doubly marked partition in order to divide the set of marked  partitions counted by $\spt(5n+4)$
(or $\spt(7n+5)$) into five (or seven) equinumerous classes. Hence, in principle, the spt-crank of a doubly marked partition can
be considered as a solution to the problem of Andrews, Dyson and Rhoades. It would be
interesting to find an spt-crank directly defined on marked partitions.

 Let $N_S(m,n)$  denote the  net number, or the sum of signs,  of $S$-partitions  of $n$ with spt-crank $m$. Andrews, Dyson and Rhoades \cite{Andrews-Dyson-Rhoades-2013} conjectured that $\{N_S(m,n)\}_m$ is unimodal for any given $n$  and showed that this conjecture is equivalent to an inequality between the rank and the crank  of a partition. Using the notion of the rank-set of a partition introduced by Dyson  \cite{Dyson-1989},    Chen, Ji and Zang \cite{Chen-Ji-Zang-2015} gave a proof of this conjecture by constructing an injection from  the set of   partitions of $n$ such that $m$ appears in the rank-set to the set of  partitions of $n$ with rank  not less than $-m$.

 \subsection{The spt-crank of an $S$-partition}

Based on \eqref{gf-spt},
 Andrews, Garvan and Liang \cite{Andrews-Garvan-Liang-2011}
 noticed that the generating
 function of $\spt(n)$ can be expressed as
 \begin{equation}\label{gf-spt-1}
 \sum_{n=1}^\infty\spt(n)q^n=\sum_{n=1}^{\infty}\frac{q^n(q^{n+1};q)_\infty}
 {(q^n;q)_\infty^2},
\end{equation}
and they introduced the structure of $S$-partitions to interpret the right-hand side of \eqref{gf-spt-1}  as the
 generating function of the net number of $S$-partitions of $n$,
 that is, the sum of   signs of $S$-partitions of $n$. More precisely, let $\mathcal{D}$ denote  the set of partitions into distinct parts and $\mathcal{P}$  denote the set of partitions. For   $\lambda \in \mathcal{P}$, we use $s(\lambda)$ to denote  the smallest part  of $\lambda$ with the convention that $s(\emptyset)=+\infty$. The set of $S$-partitions is defined by
\begin{equation}\label{def-S-par}
S=\{(\pi_1, \pi_2,\pi_3)\in \mathcal{D} \times \mathcal{P} \times \mathcal{P}\ | \  \pi_1 \neq \emptyset \text{ and }s(\pi_1)\leq \min\{s(\pi_2),s(\pi_3)\}\}.
\end{equation}
For $\pi=(\pi_1,\,\pi_2,\,\pi_3) \in S$,
 Andrews, Garvan and Liang  \cite{Andrews-Garvan-Liang-2011}  defined the weight of $\pi$ to be $|\pi_1|+|\pi_2|+|\pi_3|$ and  defined the sign of $\pi$ to be \[\omega(\pi)=(-1)^{l(\pi_1)-1},\]
 where $|\pi|$ denotes the sum of parts of $\pi$ and $l(\pi)$ denotes the number of parts of $\pi$.

 They showed that
 \[\spt(n)=\sum_{\pi }\omega(\pi),\]
 where $\pi$ ranges over $S$-partitions of $n$.
  To give combinatorial interpretations of the congruences \eqref{con-5} and \eqref{con-7},  Andrews, Garvan and Liang  \cite{Andrews-Garvan-Liang-2011} defined the   spt-crank of an $S$-partition, which  takes the same
 form as the crank of a vector partition.

 Let $\pi$ be an $S$-partition, the spt-crank of $\pi$, denoted  $r(\pi)$, is defined to be the number of parts of $\pi_2$ minus the number of parts of
  $\pi_3$, i.e.,
 \[ r(\pi)=l(\pi_2)-l(\pi_3).\]
 Let $N_S(m,n)$ denote the net number of $S$-partitions of $n$ with spt-crank $m$, that is,
\begin{equation}\label{equ-def-nsmn}
N_S(m,n)=\sum_{\stackrel{|\pi|=n}{r(\pi)=m}}\omega(\pi),
\end{equation}
and let $N_S(k,t,n)$  denote the net number of $S$-partitions of $n$ with spt-crank congruent to $k \pmod t$, namely,
\[N_S(k,t,n)=\sum_{m \equiv k \pmod{t}}N_S(m,n).\]

Andrews, Garvan and Liang \cite{Andrews-Garvan-Liang-2011}  obtained the following relations.

\begin{thm}[Andrews, Garvan and Liang]\label{AGL} For $0\leq k\leq 4$,
\begin{eqnarray*}
N_S(k,5,5n+4)&=& \frac{\spt(5n+4)}{5},
\end{eqnarray*}
and for $0\leq k\leq 6$,
\begin{eqnarray*}
N_S(k,7,7n+5)&=& \frac{\spt(7n+5)}{7}.
\end{eqnarray*}
\end{thm}

Andrews, Garvan and Liang \cite{Andrews-Garvan-Liang-2011}
defined an involution  on the set of $S$-partitions:
 \[\iota(\vec{\pi})=\iota(\pi_1,\pi_2,\pi_3)=(\pi_1,\pi_3,\pi_2),\]
which leads to the   symmetry property of $N_S(m,n)$:
\begin{equation}\label{equ-sym-nsmn}
N_S(m,n)=N_S(-m,n).
\end{equation}

Using the generating function of $N_S(m,n)$, Andrews, Garvan and Liang
\cite{Andrews-Garvan-Liang-2011}   proved its positivity.

\begin{thm}[Andrews, Garvan and Liang]For all integers $m$ and positive integers $n$,
\begin{equation}
N_S(m,n)\geq 0.
\end{equation}
\end{thm}

Dyson \cite{Dyson-2012} gave an alternative proof of this property by establishing the  relation:
\[N_S(m,n)=\sum_{k=1}^\infty (-1)^{k-1}\sum_{j=0}^{k-1}p(n-k(m+j)-(k(k+1)/2)).
\]

Andrews, Garvan and Liang  \cite{Andrews-Garvan-Liang-2011}   posed the problem of finding a combinatorial interpretation of $N_S(m,n)$.  Chen, Ji and Zang \cite{Chen-Ji-Zang-2016} introduced the structure of a doubly marked partition
 which leads to a combinatorial interpretation of $N_S(m,n)$.

 \subsection{The spt-crank of a doubly marked partition}
 In this section, we first give a definition of a doubly marked partition and then define its   spt-crank.   To this end,  we assume that a partition $\lambda$ of $n$ is represented by
its Ferrers diagram, and we use $D(\lambda)$ to denote size of the
Durfee square of $\lambda$, see \cite[p.\,28]{Andrews-1976}.   For each partition $\lambda=(\lambda_1,\lambda_2,\ldots, \lambda_l)$ of $n$, the associated Ferrers diagram
 is the arrangement of $n$ dots in $l$ rows with the dots being left-justified  and the $i$-th row having $\lambda_i$ dots for $1\leq i\leq l$. The Durfee square of $\lambda$ is the largest-size square contained within the Ferrers diagram of $\lambda$.

For a partition $\lambda$, let $\lambda'$ denote its conjugate.
A doubly marked partition of $n$ is a partition $\lambda$ of $n$ along with two distinguished columns indexed by $s$ and $t$, denoted $(\lambda,s,t)$, where
\begin{itemize}
\item[{\rm (1)}] $1\leq s\leq D(\lambda)${\rm;}
\item[{\rm (2)}] $s\leq t\leq \lambda_1${\rm;}
\item[{\rm (3)}] $\lambda'_s=\lambda'_t$.
\end{itemize}

For example, $((3,2,2),1,2)$ is a doubly marked partition, whereas \break$((3,2,1),1,2)$ and $((3,2,2),2,1)$ are not doubly marked partitions, see Figure \ref{double-1}.
\input{double.TpX}

 To define the spt-crank of a doubly marked partition  $(\lambda,s,t)$,
 let \begin{equation}\label{defi-t}
  g(\lambda,s,t)=\lambda'_{s}-s+1,
\end{equation}
  As $s\leq D(\lambda)$, we see that $\lambda'_s\geq s$, which implies that $g(\lambda,s,t)\geq 1$.

Let $(\lambda,s,t)$ be a doubly marked partition, and
 let $g=g(\lambda, s,t)$. The  spt-crank of   $(\lambda,s,t)$ is defined by
\begin{equation}\label{def-sptc}
c(\lambda,s,t)=g-\lambda_{g}+t-s.
\end{equation}

For example, for the doubly marked partition $((4,4,1,1),2,3)$, we have $g=2-1=1$ and the spt-crank  equals $1-\lambda_1+3-2=-2.$

The following theorem in \cite{Chen-Ji-Zang-2016} gives a combinatorial interpretation of $N_S(m,n)$.

\begin{thm}[Chen, Ji and Zang]\label{main-1}
For any integer  $m$ and any positive integer $n$,
 $N_S(m,n)$ equals the number of doubly marked partitions of $n$ with spt-crank $m$.
\end{thm}

 For example, for $n=4$, the sixteen $S$-partitions
of $4$, their spt-cranks  and the ten doubly marked partitions of $4$ and their spt-cranks  are listed in  Table 1.

\begin{table}[h]
\[\begin{array}{ccc|cc}
S\small \text{-partition}& \small \text{sign}&\small \text{spt-crank}&\small \text{doubly marked partition}&\small \text{spt-crank}\\[2pt]
{\footnotesize ((1),(1,1,1),\emptyset)}& \footnotesize +1&3& ((1,1,1,1),1,1)&3\\[2pt]\hline
((1),(2,1),\emptyset)&+1&2&((2,1,1),1,1)&2\\[2pt]\hline
((1),(1,1),(1))&+1&1&((3,1),1,1)&1\\[2pt]
((1),(3),\emptyset)&+1&1&((2,2),1,2)&1\\[2pt]
((2,1),(1),\emptyset)&-1&1&&\\[2pt]
((2),(2),\emptyset)&+1&1&&\\[2pt]\hline
((1),(2),(1))&+1&0&((2,2),1,1)&0\\[2pt]
((1),(1),(2))&+1&0&((4),1,4)&0\\[2pt]
((3,1),\emptyset,\emptyset)&-1&0&&\\[2pt]
((4),\emptyset,\emptyset)&+1&0&&\\[2pt]\hline
((1),(1),(1,1))&+1&-1&((2,2),2,2)&-1\\[2pt]
((1),\emptyset,(3))&+1&-1&((4),1,3)&-1\\[2pt]
((2,1),\emptyset,(1))&-1&-1&&\\[2pt]
((2),\emptyset,(2))&+1&-1&&\\[2pt]\hline
((1),\emptyset,(2,1))&+1&-2&((4),1,2)&-2\\[2pt]\hline
((1),\emptyset,(1,1,1))&+1&-3&((4),1,1)&-3
\end{array}
\]
\caption{$S$-partitions and doubly marked partitions}
\end{table}

The proof of  Theorem \ref{main-1} relies on the generating function of $N_S(m,n)$ given by Andrews, Garvan and Liang \cite{Andrews-Garvan-Liang-2011}.

  Andrews, Dyson and Rhoades \cite{Andrews-Dyson-Rhoades-2013}  proposed the problem of  finding a definition of the spt-crank for a marked partition so that the set of marked partitions of $5n+4$ (or $7n+5$) can be divided into five (or seven) equinumerous classes.
 Chen, Ji and Zang \cite{Chen-Ji-Zang-2016} established a bijection $\Delta$ between the set of marked partitions of $n$ and the set of doubly marked partitions of $n$.

\begin{thm}[Chen, Ji and Zang]\label{main-2} There is a bijection $\Delta$ between the set of marked partitions $(\mu,k )$ of $n$ and the set of doubly marked partitions $(\lambda, s, t)$ of $n$.
\end{thm}

To prove the above theorem, we adopt the notation $(\lambda, s,t)$
 for a partition $\lambda$ with two distinguished columns $\lambda'_s$ and $\lambda'_t$ in the Ferrers diagram.
Let $Q_n$ denote the set of doubly marked partitions of $n$,
 and let
 \[U_n=\{(\lambda,s,t) \ |\    |\lambda|=n,\, 1\leq s\leq D(\lambda),\, 1\leq t\leq \lambda_1\}.\]
 Obviously, $Q_n\subseteq U_n$.

Before we give a description of the bijection $\Delta$, we  introduce a transformation $\tau$ from   $U_n\setminus Q_n$ to $U_n$.

\noindent{\it The transformation $\tau$:} Assume that $(\lambda,s,t) \in U_n\setminus Q_n$, that is, $\lambda$ is a partition of $n$ with
two distinguished columns indexed by $s$ and $t$ such that $1\leq s\leq D(\lambda)$  and either $1\leq t<s$ or $\lambda'_{s}>\lambda'_t.$  We wish to construct a partition $\mu$ with two distinguished columns  indexed by $a$ and $b$.  Let  $p$ be the maximum integer such that $\lambda_p'=\lambda_s'$. Define
\begin{equation}\label{phi-gamma}
\delta=(\lambda_1-p+s-1,\lambda_2-p+s-1,\ldots, \lambda_{\lambda_s'}-p+s-1,\lambda_{\lambda_s'+1},\ldots, \lambda_\ell).
\end{equation}
 Set $a$ to be   the minimum integer such that $\delta_a<\lambda'_s$ and
\begin{equation}\label{def-mu}
\mu=(\delta_1,\ldots,\delta_{a-1},\lambda'_s,\ldots, \lambda'_p, \delta_{a},\ldots, \delta_\ell).
\end{equation}
If $t<s$, then set $b=t$ and if $\lambda_s'>\lambda_t'$, then   set $b=t-p+s-1$.  Define $\tau(\lambda,s,t)=(\mu,a,b)$. Figure \ref{fig-phi} gives an illustration of the map $\tau\colon((6,5,3,1),2,6)\mapsto ((4,3,3,3,1,1),3,4)$.
  \input{phi.TpX}

It was proved in \cite{Chen-Ji-Zang-2016}  that the map $\tau$ is indeed an injection. Using this property, they described the bijection $\Delta$ in Theorem \ref{main-2} based on the injection $\tau$.

\noindent{\it The definition of $\Delta \colon$} Let $(\mu, k)$ be a marked partition of $n$, we proceed to construct a doubly marked partition $(\lambda,s,t)$ of $n$.

 We first consider  $(\mu', 1, k)$.
 If $(\mu', 1, k)$ is already a doubly marked partition, then
 there is nothing to be done and we just set $(\lambda,s,t)=(\mu',1,k)$. Otherwise, we iteratively apply the map $\tau$ to $(\mu',1,k)$
 until we get a doubly marked partition $(\lambda,s,t)$. We then define
 \[\Delta(\mu,k)=(\lambda,s,t).\]
 It can be shown that
 this process terminates and it is reversible. Thus $\Delta$ is well-defined and is  a bijection between the set of marked partitions $(\mu,k )$ of $n$ and the set of doubly marked partitions $(\lambda, s, t)$ of $n$.

 To give an example of the map $\Delta$, let $n=6$, $\mu=(2,1,1,1,1)$ and $k=5$. We have $\mu'=(5,1)$.
 Note that $(\mu', 1, k)=((5,1),1,5)$, which  is not a doubly marked partition.
 It can be checked that  $\tau(\mu', 1,k)=
       ((4,2),2,4)$, which is not a doubly marked partition.
 Repeating this process, we get $\tau((4,2),2,4)=((3,2,1),2,3)$,
  and $\tau ((3,2,1),2,3)=((2,2,1,1),2,2)$, \break which is eventually  a doubly marked partition.  See Figure \ref{fig2}. Thus, we obtain
 \[\Delta((2,1,1,1,1),5)=((2,2,1,1),2,2).\]

{\small
 \input{triexa.TpX}
}

Utilizing the bijection $\Delta$ and the spt-crank for a doubly marked partition, one can
divide the set of marked partitions of $5n + 4$ (or $7n + 5$) into five (or seven) equinumerous classes. Hence, in principle,
the spt-crank of a doubly marked partition can
be considered as a solution to the problem of Andrews, Dyson and Rhoades.
It would be interesting to find an spt-crank  directly
defined on marked partitions.

For example, for $n=4$, we have $\spt(4)=10$. The ten marked partitions
of $4$, the corresponding doubly marked partitions, and the spt-crank modulo $5$ are listed in Table \ref{table-1}.

{\small
 \begin{table}[!h]
 \[\begin{array}{c|c|c|c}
(\mu,k)& (\lambda,s,t)=\Delta(\mu,k) & c(\lambda,s,t) &c(\lambda,s,t) \mod 5\\[2pt] \hline
((4),1) & ((1,1,1,1),1,1)& 3& 3\\[2pt]
((3,1),2)& ((3,1),1,1)& 1& 1\\[2pt]
((2,2),1)& ((2,2),1,1)& 0& 0\\[2pt]
((2,2),2)& ((2,2),1,2)& 1& 1\\[2pt]
 ((2,1,1),2)& ((2,1,1),1,1)& 2 &2\\[2pt]
((2,1,1),3)& ((2,2),2,2)& -1 &4\\[2pt]
((1,1,1,1),1)& ((4),1,1) & -3 &2\\[2pt]
((1,1,1,1),2)& ((4),1,2) & -2 &3\\[2pt]
((1,1,1,1),3)& ((4),1,3) & -1 &4\\[2pt]
((1,1,1,1),4)& ((4),1,4) & 0 &0
\end{array}\]

 \caption{ The case for $n=4$} \label{table-1}
\end{table}
}

For $n=5$, we have $\spt(5)=14$. The fourteen marked partitions
of $5$, the corresponding doubly marked partitions, and the spt-crank modulo $7$ are listed in  Table \ref{table-2}.

{\small
\begin{table}[!h]
\[\begin{array}{c|c|c|c}
(\mu,k)& (\lambda,s,t)=\Delta(\mu,k) & c(\lambda,s,t) & c(\lambda,s,t)\mod 7\\[2pt] \hline
((5),1) & ((1,1,1,1,1),1,1)& 4& 4\\[2pt]
((4,1),2)& ((4,1),1,1)& 1& 1\\[2pt]
((3,2),2)& ((3,1,1),1,1)& 2& 2\\[2pt]
((3,1,1),2)& ((3,2),1,1)& 0& 0\\[2pt]
((3,1,1),3)& ((3,2),1,2)& 1& 1\\[2pt]
((2,2,1),3)& ((2,2,1),1,1)& 2& 2\\[2pt]
((2,1,1,1),2)& ((2,1,1,1),1,1)& 3 &3\\[2pt]
((2,1,1,1),3)& ((3,2),2,2)& -2 &5\\[2pt]
((2,1,1,1),4)& ((2,2,1),2,2)& -1 &6\\[2pt]
((1,1,1,1,1),1)& ((5),1,1) & -4 &3\\[2pt]
((1,1,1,1,1),2)& ((5),1,2) & -3 &4\\[2pt]
((1,1,1,1,1),3)& ((5),1,3) & -2 &5\\[2pt]
((1,1,1,1,1),4)& ((5),1,4) & -1 &6\\[2pt]
((1,1,1,1,1),5)& ((5),1,5) & 0 &0
\end{array}\]

 \caption{ The case for $n=5$} \label{table-2}
\end{table}
}

\subsection{The unimodality of  the spt-crank}

  The unimodality of the spt-crank was first studied by Andrews, Dyson and Rhoades \cite{Andrews-Dyson-Rhoades-2013}. They showed that the unimodality of the spt-crank is equivalent to an inequality between  the rank and the crank of a partition. Define
\begin{eqnarray}
N_{\leq m}(n)&=&\sum_{|r|\leq m}N(r,n), \\[3pt]
M_{\leq m}(n)&=&\sum_{|r|\leq m}M(r,n).
\end{eqnarray}

Andrews, Dyson and Rhoades \cite{Andrews-Dyson-Rhoades-2013} established the following relation.

\begin{thm}[Andrews, Dyson and Rhoades]\label{relation-spt-crank}
For  $m\geq 0$  and \break $n>1$,
 \begin{equation}\label{equiv-e}
 N_S(m,n)- N_S(m+1,n)=\frac{1}{2}\left(N_{\leq m}(n)- M_{\leq m}(n)\right).
 \end{equation}
 \end{thm}

They also  posed a conjecture on the spt-crank.

   \begin{conje}[Andrews, Dyson and Rhoades]\label{conj-o} For $m,n\geq 0$,
 \begin{equation}\label{ine-conj-nsum}
 N_S(m,n)\geq N_S(m+1,n).
 \end{equation}
 \end{conje}

By the symmetry   \eqref{equ-sym-nsmn} of $N_S(m,n)$ and the
relation \eqref{ine-conj-nsum}, we see that
\[N_S(-n,n)\leq\cdots\leq N_S(-1,n)\leq N_S(0,n)\geq N_S(1,n)\geq\cdots\geq N_S(n,n).\]

{\tiny
 \begin{table}[h!]
\[\begin{array}{c|ccccccccccccccccccccccccccccc}
n\setminus m&-6&-5&-4&-3&-2&-1&\ 0&\ 1&\ 2&\ 3&\ 4&\ 5&\ 6\\[2pt]\hline
0&&&&&&&1&&&&&&\\[2pt]
1&&&&&&&1&&&&&&\\[2pt]
2&&&&&&1&1&1&&&&&\\[2pt]
3&&&&&1&1&1&1&1&&&&\\[2pt]
4&&&&1&1&2&2&2&1&1&&&\\[2pt]
5&&&1&1&2&2&2&2&2&1&1&&\\[2pt]
6&&1&1&2&3&4&4&4&3&2&1&1&&&\\[2pt]
7&1&1&2&3&4&4&5&4&4&3&2&1&1&
\end{array}
\]
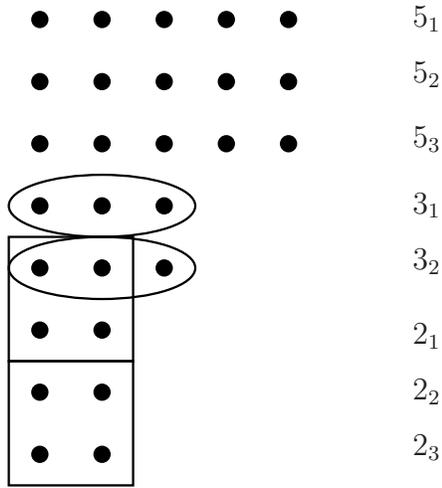
 \caption{ An illustration of the unimodality of $N_S(m,n)$}
\end{table}
}
In view of (\ref{equiv-e}), Andrews, Dyson and Rhoades pointed out that Conjecture \ref{conj-o} is equivalent to
the assertion
\begin{equation}\label{conj-e}
N_{\leq m}(n)\geq M_{\leq m}(n),
\end{equation}
where  $m, n\geq0$. It was  remarked in \cite{Andrews-Dyson-Rhoades-2013} that \eqref{conj-e} was   conjectured by Bringmann and Mahlburg \cite{Bringmann-Mahlburg-2009}. When $m=0$,   (\ref{conj-e}) was  conjectured by Kaavya \cite{Kaavya-2011}.

 Andrews, Dyson and Rhoades \cite{Andrews-Dyson-Rhoades-2013}  obtained  an asymptotic formula for $N_{\leq m}(n)-M_{\leq m}(n)$, which implies that  Conjecture \ref{conj-o} holds for fixed $m$ and sufficiently large $n$.

\begin{thm}[Andrews, Dyson and Rhoades] For each $m\geq 0$,
\begin{equation}
(N_{\leq m}(n)-M_{\leq m}(n))\sim \frac{(2m+1)\pi^2}{192\sqrt{3}n^2}e^{\pi\sqrt{\frac{2n}{3}}} \quad \text{as} \quad n\rightarrow \infty.
\end{equation}
\end{thm}

 Using the rank-set of a partition, Chen, Ji and Zang   \cite{Chen-Ji-Zang-2015} constructed  an injection from  the set of  partitions of $n$ such that $m$ appears in the rank-set to the set of  partitions of $n$ with rank  not less than $-m$. This proves the inequality \eqref{conj-e}   for all $m\geq 0$ and $n\geq 1$, and hence Conjecture \ref{conj-o} is confirmed.

In fact, the relation \eqref{conj-e} was stated by Bringmann and Mahlburg \cite{Bringmann-Mahlburg-2009}
in a different notation.  For an integer $m$ and a positive integer $n$, let
 \[\overline{\mathcal{M}}(m,n) =\sum_{r\leq m}M(r,n),\]
 and
  \[\overline{\mathcal{N}}(m,n) =\sum_{r\leq m}N(r,n).\]
By the symmetry properties of the rank and the crank, that is,
\[N(m,n)=N(-m,n)\quad\text{and}\quad M(m,n)=M(-m,n),\]
see \cite{Dyson-1944} and \cite{Garvan-1988},
it is not difficult to verify that \eqref{conj-e} is equivalent to the following inequality for $m<0$ and $n\geq 1$:
 \begin{equation}\label{conj-bm-e}
 \overline{\mathcal{N}}(m,n)\leq \overline{\mathcal{M}}(m,n).
 \end{equation}

 It turns out that the constructive approach in \cite{Chen-Ji-Zang-2015} can be used to prove the other part of the conjecture \eqref{conj-bm-e} of
Bringmann and Mahlburg, that is,
\begin{equation}\label{BM-CI}
 \overline{\mathcal{M}}(m,n)\leq \overline{\mathcal{N}}(m+1,n),
 \end{equation}
 for $m< 0$ and $n\geq 1$.
 A  proof of  \eqref{BM-CI} was given in \cite{Chen-Ji-Zang-2016p}.

In the notation    $N_{\leq m-1}(n)$ and $M_{\leq m}(n)$, the inequality \eqref{BM-CI}  can be expressed as
\begin{equation}\label{conj-BM-eqn}
  M_{\leq m}(n)\geq N_{\leq m-1}(n),
  \end{equation}
  for $m\geq 1$ and $n\geq 1$.

 Bringmann and Mahlburg \cite{Bringmann-Mahlburg-2009}  also pointed out that the inequalities \eqref{conj-bm-e} and \eqref{BM-CI} can be restated as the existence of a re-ordering $\tau_n$ on the set of partitions of $n$ such that $|\text{crank}(\lambda)|-|\text{rank}(\tau_n(\lambda))|=0$ or $1$ for all partitions $\lambda$ of $n$.  Chen, Ji and Zang \cite{Chen-Ji-Zang-2016p}   defined a  re-ordering $\tau_n$
  on the set of partitions of $n$ and showed that this re-ordering  $\tau_n$ satisfies the relation $|\text{crank}(\lambda)|-|\text{rank}(\tau_n(\lambda))|=0$ or $1$ for any partition
  $\lambda$ of $n$. Appealing to this re-ordering $\tau_n$,  they  gave a new combinatorial interpretation of  the
function $\ospt(n)$ defined by Andrews, Chan and Kim \cite{Andrews-Chan-Kim-2013}, which  leads to  an upper bound for  $\ospt(n)$ due to  Chan and Mao \cite{Chan-Mao-2014}.

 Bringmann and Mahlburg \cite{Bringmann-Mahlburg-2009} also remarked that using the Cauchy-Schwartz inequality,  the bijection $\tau_n$  leads to an upper bound for $\spt(n)$, namely, for $n\geq 1$,
 \begin{equation}\label{ine-spt-2npn}
 \spt(n)\leq \sqrt{2n}p(n).
 \end{equation}

Chan and Mao \cite{Chan-Mao-2014} posed a conjecture on a  sharper upper bound and a lower bound  for $\spt(n)$.
\begin{conj}[Chan and Mao]
For $n\geq 3$,
\begin{equation}
\frac{\sqrt{6n}}{\pi} p(n)\leq \spt(n)\leq \sqrt{n}p(n).
\end{equation}
\end{conj}

The following upper bound and lower bound for $\spt(n)$ were conjectured by Hirschhorn  and later proved by Eichhorn and Hirschhorn  \cite{Eichhorn-Hirschhorn-2015}.

\begin{thm}[Eichhorn and Hirschhorn]
For $n\geq 2$,
\begin{equation}
p(0)+p(1)+\cdots+p(n-1)< \spt(n)< p(0)+p(1)+\cdots+p(n).
\end{equation}
\end{thm}

\section{More congruences}

Garvan  \cite{Garvan-2010} obtained     Ramanujan-type
congruences of $\spt(n)$ mod $11$, $17$, $19$, $29$, $31$ and $37$.

\begin{thm}[Garvan]
For $n\geq 0$,
\begin{eqnarray}
\spt(11\cdot19^4\cdot n+22006)&\equiv& 0\pmod{11},\\[3pt]
 \spt(17\cdot7^4\cdot n+243)&\equiv& 0\pmod{17},\\[3pt]
 \spt(19\cdot5^4\cdot n+99)&\equiv& 0\pmod{19},\\[3pt]
 \spt(29\cdot13^4\cdot n+18583)&\equiv& 0\pmod{29},\\[3pt]
 \spt(31\cdot29^4\cdot n+409532)&\equiv& 0\pmod{31},\\[3pt]
 \spt(37\cdot5^4\cdot n+1349)&\equiv& 0\pmod{37}.
\end{eqnarray}
\end{thm}

Bringmann \cite{Bringmann-2008} showed that $\spt(n)$
possesses a congruence property analogous to the following theorem   for $p(n)$, due to Ono \cite{Ono-2000}.

  \begin{thm}[Ono]\label{thm-o-pn}
 For any prime $\ell\geq 5$,  there are infinitely  many arithmetic progressions $an+b$ such that
  \begin{equation}\label{equ-mod-pn-ar}
  p(an+b)\equiv 0\pmod{\ell}.
  \end{equation}
  \end{thm}

As for $\spt(n)$, Bringmann \cite{Bringmann-2008} proved the following assertion.

\begin{thm}[Bringmann]
For any prime $\ell\geq 5$, there are infinitely many
arithmetic progressions $an+b$ such that
\[\spt(an+b)\equiv 0\pmod{\ell}.\]
\end{thm}

 The above theorem is a consequence of \eqref{spt-moments},  Theorem \ref{thm-o-pn} and the following theorem of Bringmann \cite{Bringmann-2008}.

 \begin{thm}[Bringmann]\label{thm-b-n2}
 For any prime $\ell\geq 5$,  there are infinitely  many arithmetic progressions $an+b$ such that
 \begin{equation}\label{equ-mod-n2n}
 N_2(an+b)\equiv 0\pmod{\ell}.
 \end{equation}
 \end{thm}

  Bringmann \cite{Bringmann-2008} constructed a weight 3/2 harmonic weak Maass form $\mathcal{M}(z)$ on $\Gamma_0(576)$ with Nebentypus $\chi_{12}(\bullet)=\left(\frac{12}{\bullet}\right)$, which is related to the generating function of $\spt(n)$. This implies that the generating function of $\spt(n)$ is essentially a mock theta function with Dedekind eta-function $\eta(q)$ as its shadow just as pointed out by Rhoades \cite{Rhoades-2013-1}.  Ono \cite{Ono-2011} found a weight $(\ell^2 + 3)/2$ holomorphic modular form on $SL_2(\mathbb{Z})$ which contains the holomorphic part of $\mathcal{M}(z)$. Using this modular form,
Ono \cite{Ono-2011} derived  Ramanujan-type congruences of $\spt(n)$ modulo $\ell$ for any prime $\ell\geq 5$.

\begin{thm}[Ono]\label{Ono-con}
Let $\ell\geq 5$ be a prime and let $\left(\frac{\bullet}{\circ}\right)$ denote the Legendre symbol.
 \begin{itemize}
 \item[{\rm(i)}] For $n\geq 1$, if $\left(\frac{-n}{\ell}\right)=1$,
\[\spt\left((\ell^2 n+1)/24\right)\equiv 0\pmod{\ell}.
\]
\item[{\rm(ii)}] For $n\geq 0$,
\[\spt\left((\ell^{3}n+1)/24\right)\equiv \left(\frac{3}{\ell}\right)\spt\left((\ell n+1)/24\right)\pmod{\ell}.\]
\end{itemize}
\end{thm}

Ahlgren,  Bringmann and  Lovejoy \cite{Ahlgren-Bringmann-Lovejoy-2011} extended  Theorem \ref{Ono-con} to any prime power. An analogous  congruence     for $p(n)$  was found by Ahlgren \cite{Ahlgren-2000}.

\begin{thm}[Ahlgren, Bringmann and Lovejoy]\label{thm-ell-adic}
Let $\ell\geq 5$ be a prime and let $m\geq 1$.
\begin{itemize}
\item[{\rm(i)}] For $n\geq 1$, if $\left(\frac{-n}{\ell}\right)=1$,
\[\spt\left((\ell^{2m}n+1)/24\right)\equiv 0\pmod{\ell^m}.\]

\item[{\rm(ii)}] For $n\geq 0$,
\[\spt\left((\ell^{2m+1}n+1)/24\right)\equiv\left(\frac{3}{\ell}\right) \spt\left((\ell^{2m-1}n+1)/24\right)\pmod{\ell^m}.\]
\end{itemize}
\end{thm}

Recall the following congruences of $p(n)$:
\begin{eqnarray}
p(5^an+\delta_a)&\equiv& 0\pmod{5^a},\label{equ-mod-5-p}\\[3pt]
p(7^bn+\lambda_b)&\equiv& 0\pmod{7^{\left\lfloor\frac{b+2}{2}\right\rfloor}},\label{equ-mod-7-p}\\[3pt]
p(11^cn+\varphi_c)&\equiv& 0\pmod{11^c},\label{equ-mod-11-p}
\end{eqnarray}
 where $a,b,c$ are positive integers and $\delta_a,\lambda_b$ and $\varphi_c$ are the least nonnegative residues of the reciprocals of $24$ mod $5^a,7^b$ and $11^c$, respectively. The congruences  \eqref{equ-mod-5-p} and \eqref{equ-mod-7-p} were proved by Watson \cite{Watson-1938} and    the  congruence \eqref{equ-mod-11-p} was proved by Atkin \cite{Atkin-1967}. Folsom, Kent and Ono \cite{Folsom-Kent-Ono-2012} provided  alternative proofs of the congruences  \eqref{equ-mod-5-p}--\eqref{equ-mod-11-p} with the aid of   the theory of $\ell$-adic modular forms. Recently, Paule and Radu \cite{Paule-Radu-2017} found a unified algorithmic
 approach to   \eqref{equ-mod-5-p}--\eqref{equ-mod-11-p}  resorting to  elementary modular function tools only.

In the case of the spt-function, although  Theorem \ref{thm-ell-adic} gives congruences for all primes $\ell\geq 5$, the congruences \eqref{con-5}--\eqref{con-13} do not follow
   from  Theorem \ref{thm-ell-adic}. Congruences for these missing cases
 have been obtained by  Garvan  \cite{Garvan-2012}, which are analogous to
 \eqref{equ-mod-5-p}--\eqref{equ-mod-11-p}.

\begin{thm}[Garvan]\label{thm-Garvan-5713p}
For $n\geq0$,
\[\spt(5^an+\delta_a)\equiv 0\pmod{5^{\lfloor\frac{a+1}{2}\rfloor}},\]
\[\spt(7^bn+\lambda_b)\equiv 0\pmod{7^{\lfloor\frac{b+1}{2}\rfloor}},\]
\[\spt(13^cn+\gamma_c)\equiv 0\pmod{13^{\lfloor\frac{c+1}{2}\rfloor}},\]
where $a,b,c$ are positive integers, and $\delta_a$, $\lambda_b$ and $\gamma_c$ are the least nonnegative residues of the reciprocals of $24$ mod $5^a$, $7^b$ and $13^c$ respectively.
\end{thm}

Setting $a=b=c=1$, Theorem \ref{thm-Garvan-5713p} reduces to  \eqref{con-5}--\eqref{con-13}. Belmont, Lee,  Musat and Trebat-Leder \cite{Belmont-Lee-Musat-2014} provided another proof of the above theorem by  generalizing   techniques of Folsom, Kent and Ono \cite{Folsom-Kent-Ono-2012} and by utilizing  refinements due to Boylan and Webb \cite{Boylan-Webb-2013}.

Before we get into the discussions about the parity of $\spt(n)$,
let us look back at the parity of $p(n)$. Subbarao \cite{Subbarao-1966} conjectured
that in every arithmetic progression $r \pmod t$,  there are infinitely many
integers $N\equiv r \pmod t$ for which $p(N)$ is even, and infinitely many integers
$M\equiv r \pmod t$ for which $p(M)$ is odd. This conjecture has been confirmed
   for  $t=1,2,3,4,5,10, 12,16$ and $40$ by Garvan and Stanton \cite{Garvan-Stanton-1990}, Hirschhorn \cite{Hirschhorn-1993}, Hirschhorn and Subbarao \cite{Hirschhorn-Subbarao-1988}, Kolberg \cite{Kolberg-1959} and Subbarao \cite{Subbarao-1966}. The even case of Subbarao's conjecture was
settled by Ono \cite{Ono-1996} and  the odd case
 was solved by Radu \cite{Radu-2012}. Radu \cite{Radu-2012} also  showed
that for every arithmetic progression $r \pmod t$, there are infinitely many integers
$N\equiv r \pmod t$ such that $p(N)\not\equiv 0 \pmod 3$. This confirms a
conjecture posed by Ahlgren and Ono \cite{Ahlgren-Ono-2001}.

For $n\geq 1$, the parity of $\spt(n)$ is determined by Folsom and Ono \cite{Folsom-Ono-2008}. They
constructed a pair of harmonic weak Maass forms with equal  nonholomorphic parts, whose difference contains the generating function of $\spt(n)$ as a component. Based on the results in \cite{Bringmann-Folsom-Ono-2009}, Folsom and Ono showed that the difference of such pair of harmonic weak Maass forms can be expressed as
the sum of the generating function for $\spt(n)$ and a modular form. This enables us to completely determine  the parity of $\spt(n)$.

To be more specific, Folsom and Ono \cite{Folsom-Ono-2008} first defined the  mock theta functions:
\[D(z) =\frac{q^{-\frac{1}{24}}}
{(q;q)_\infty}\left(1-24\sum_{n=1}^\infty
\frac{nq^n}{1-q^n}\right)=\frac{q^{-\frac{1}
{24}}}{(q;q)_\infty}E_2(z)\]
and
{\small
\[L(z) =\frac{(q^6;q)_\infty^2(q^{24};q)_\infty^2}
{(q^{12};q)^5_\infty}\left(
\sum_{n=-\infty}^{\infty}\frac{(12n-1)q^{6n^2-\frac{1}{24}}}{1-q^{12n-1}}
-\sum_{n=-\infty}^{\infty}\frac{(12n-5)q^{6n^2-\frac{25}{24}}}{1-q^{12n-5}}
\right).\]}

Then they  obtained the following modular form.

\begin{thm}[Folsom and Ono]\label{thm-modu-spt}
The function
\[D(24z)-12L(24z)-12q^{-1}S(24z)\]
is a weight $3/2$ weakly holomorphic modular form on $\Gamma_0(576)$ with Nebentypus $\left(\frac{12}{\bullet}\right)$, where
\[S(z)=\sum_{n=0}^\infty \spt(n)q^n.\]
\end{thm}

By Theorem \ref{thm-modu-spt}, Folsom and Ono \cite{Folsom-Ono-2008}
obtained a characterization of the  parity of $\spt(n)$.

\begin{thm}[Folsom and Ono]\label{thm-spt-parity}
 The function $\spt(n)$ is odd if and only if $24n-1=pm^2$, where $m$ is an integer and $p\equiv 23\pmod{24}$ is  prime.
\end{thm}

As pointed out by Andrews, Garvan and Liang  \cite{Andrews-Garvan-Liang-2013},
 Theorem \ref{thm-spt-parity} contains an error. For example, for $n=507$, it is clear that $507\times 24-1=12167=23\times 23^2=pm^2$, where $p=m=23$. Obviously, $507$ satisfies the condition of Theorem \ref{thm-spt-parity}. But  $\spt(507)=60470327737556285225064$ is even. This error has been corrected by
  Andrews, Garvan and Liang \cite{Andrews-Garvan-Liang-2013}.
   By using the notion of  $S$-partitions as defined in \eqref{def-S-par},
  they noticed that the number of $S$-partitions of $n$  has the same parity as $\spt(n)$. Then they built an involution $\iota$ on the set of $S$-partitions of $n$ as follows:
 \[\iota(\vec{\pi})=\iota(\pi_1,\pi_2,\pi_3)=(\pi_1,\pi_3,\pi_2).\]
 Clearly, an $S$-partition $(\pi_1,\pi_2,\pi_3)$ is a fixed point of $\iota$ if and only if $\pi_2=\pi_3$. Denote the number of such $S$-partitions of $n$ by $N_{SC}(n)$. It is not difficult to see that
 \[\spt(n)\equiv N_{SC}(n)\pmod{2}.\]

 By computing the generating function of $N_{SC}(n)$,  Andrews, Garvan and Liang \cite{Andrews-Garvan-Liang-2013} established  a  corrected version of Theorem \ref{thm-spt-parity}.

\begin{thm}[Andrews, Garvan and Liang]\label{thm-pir-spt-co}
The function $\spt(n)$ is odd if and only if $24n-1 = p^{4a+1}m^2$ for some prime $p \equiv 23  \pmod{24}$
and some integers $a, m$ with $(p,m) = 1$.
\end{thm}

The spt-function is also related to some combinatorial sequences, see, for example,  Andrews, Rhoades and Zwegers \cite{Andrews-Rhoades-Zwegers-2013} and Bryson, Ono, Pitman and Rhoades \cite{Bryson-Ono-Pitman-Rhoades-2012}.   Bryson, Ono, Pitman and Rhoades \cite{Bryson-Ono-Pitman-Rhoades-2012} showed that the number of strongly unimodal sequences of size $n$ has the same parity as $\spt(n)$. More specifically, a sequence of integers $\{a_i\}_{i=1}^s$ is
said to be a strongly unimodal sequence of size $n$ if $a_1+\cdots+a_s=n$ and for some $k$,
\[0<a_1<a_2<\cdots<a_k>a_{k+1}>a_{k+2}>\cdots>a_s>0.\]
Let $u(n)$ be the number of strongly unimodal sequences of size $n$. By \cite[Theorem 1]{Andrews-2013}, Bryson, Ono, Pitman and Rhoades \cite{Bryson-Ono-Pitman-Rhoades-2012} observed that
\[u(n) \equiv \spt(n) \pmod{2}. \]

As for congruences of $\spt(n)$ modulo powers of $2$,
 Garvan and Jennings-Shaffer \cite{Garvan-Jennings-2014-1}
   obtained congruences mod $2^3, 2^4$ and $2^5$. Let
   $$s_\ell =\frac{\ell^2-1}{24}.$$
  % and
%   \begin{equation*}
%   \chi_{12}(\ell)=\begin{cases}
%   1&\text{if }\ell\equiv\pm 1\pmod{12};\\
%   -1&\text{if }\ell\equiv\pm 5\pmod{12};\\
%   0&\text{otherwise.}
%   \end{cases}
%   \end{equation*}

\begin{thm}[Garvan and Jennings-Shaffer]
Let $\ell\geq 5$ be a prime, and define
\[\beta=\begin{cases}
3,&\text{if }\ell\equiv 7,9\pmod{24},\\[3pt]
4,&\text{if }\ell\equiv 13,23\pmod{24},\\[3pt]
5,&\text{if }\ell\equiv 1,11,17,19\pmod{24}.
\end{cases}\]
Then for $n\geq 1$,
\begin{gather*}
\spt(\ell^2n-s_\ell)+ \left(\frac{3-72n}{\ell}\right)\spt(n)+\ell \spt\left((n+s_\ell)/\ell^2\right)\\[3pt]
\equiv \left(\frac{3}{\ell}\right) (1+\ell)\spt(n)\pmod{2^\beta}.
\end{gather*}
\end{thm}

By using the Hecke algebra of a Maass
form,
Folsom and Ono \cite{Folsom-Ono-2008} derived a
congruence of $\spt(n)$ modulo $3$.

\begin{thm}[Folsom and Ono]\label{thm-mod-3}
Let $\ell\geq 5$ be a prime, then for $n\geq1$,
\begin{gather*}
\spt(\ell^2n-s_\ell)+\left(\frac{3-72n}{\ell}\right)\spt(n)+\ell \spt\left((n+s_\ell)/\ell^2\right)\\[3pt]
\equiv
\left(\frac{3}{\ell}\right) (1+\ell)\spt(n)\pmod{3}.
\end{gather*}
\end{thm}

\begin{core}[Folsom and Ono]
Let $\ell\geq 5$ be a prime such that $\ell\equiv 2\pmod{3}$. If $0 < k < \ell-1,$
then for   $n\geq 1$,
\begin{eqnarray*}
&&\spt(\ell^4n+\ell^3k-(\ell^4-1)/24)\equiv
 0 \pmod{3}.\end{eqnarray*}
\end{core}

For example, for $\ell=5$, we have
\begin{eqnarray*}
\spt(625n + 99) &\equiv& \spt(625n + 224)\\[3pt]
&\equiv& \spt(625n + 349)\\[3pt]
 &\equiv& \spt(625n + 474)\\[3pt]
 & \equiv& 0 \pmod {3}.\end{eqnarray*}

Garvan \cite{Garvan-2013} derived congruences mod $5, 7, 13$ and  $72$.

\begin{thm}[Garvan]
{\rm (i)} If $\ell\geq 5$ is prime, then for $n\geq1$
\begin{gather}
\spt(\ell^2n-s_\ell)+\left(\frac{3-72n}{\ell}\right)\spt(n)+\ell \spt\left(({n+s_\ell})/\ell^2\right)\nonumber\\[3pt]
 \equiv
\left(\frac{3}{\ell}\right) (1+\ell)\spt(n)\pmod{72}.\label{equ-mod-spt-72}
\end{gather}

{\rm (ii)} If $\ell\geq 5$ is prime, $t=5,7$ or $13$ and $\ell\neq t$, then for $n\geq1$
\begin{gather}
\spt(\ell^2n-s_\ell)+\left(\frac{3-72n}{\ell}\right)\spt(n)+\ell \spt\left((n+s_\ell)/\ell^2\right)\nonumber\\[3pt]
 \equiv
\left(\frac{3}{\ell}\right) (1+\ell)\spt(n)\pmod{t}.\label{equ-mod-spt-5713}
\end{gather}
\end{thm}

Note that Theorem \ref{thm-mod-3} can be deduced from \eqref{equ-mod-spt-72}. Moreover, writing $32760=2^3\cdot3^2\cdot5\cdot7\cdot13$, from \eqref{equ-mod-spt-72} and \eqref{equ-mod-spt-5713},
 it is easy to deduce a congruence of $\spt(n)$ modulo $32760$.

\begin{core}[Garvan]
If $\ell$ is prime and $\ell\not\in\{ 2,3,5,7,13\}$, then for $n\geq1$
\begin{gather*}
\spt(\ell^2n-s_\ell)+\left(\frac{3-72n}{\ell}\right)\spt(n)+\ell \spt\left((n+s_\ell)/\ell^2\right)\nonumber\\[3pt]
 \equiv
\left(\frac{3}{\ell}\right)(1+\ell)\spt(n)\pmod{32760}.
\end{gather*}
\end{core}

 Garrett, McEachern, Frederick and Hall-Holt \cite{Garrett-McEachern-Frederick-2011}
  obtained a recurrence relation for $\spt(n)$. To compute $\spt(n)$, they
 introduced two integer arrays $A(n,j)$ and $B(n,j)$, where $A(n,j)$ denotes the number of  partitions of $n$ with the smallest part at least $j$ and $B(n,j)$ denotes the number of times that $j$ occurs as the smallest part of partitions of $n$. From the definitions of $A(n,j)$ and $B(n,j)$, it is not difficult to deduce the following recurrence relations:
 \begin{eqnarray*}
 A(n,j)&=& A(n-j,j)+A(n,j+1),\\[3pt]
 B(n,j)&=& A(n-j,j)+B(n-j,j),
 \end{eqnarray*}
  where $A(n,j)=B(n,j)=0$ whenever $n<j$ and $A(n,n)=B(n,n)=1$.

Thus we have
 \[\spt(n)=\sum_{j=1}^nB(n,j).\]
By the above relation, Garrett, McEachern, Frederick and Hall-Holt computed  the first
million values of $\spt(n)$, and found many conjectures on congruences
of $\spt(n)$.
\begin{eqnarray}\label{Ono-pf-cong-1}
\spt(1331n+479)&\equiv&0\pmod{11},\\[3pt]
\spt(1331n+842)&\equiv&0\pmod{11},\\[3pt]
\spt(1331n+1084)&\equiv&0\pmod{11},\\[3pt]
\spt(1331n+1205)&\equiv&0\pmod{11},\\[3pt]
\spt(1331n+1326)&\equiv&0\pmod{11},\\[3pt]
\spt(4913n+566)&\equiv&0\pmod{17},\\[3pt]
\spt(4913n+2300)&\equiv&0\pmod{17},\\[3pt]
\spt(4913n+2878)&\equiv&0\pmod{17},\\[3pt]
\spt(4913n+3167)&\equiv&0\pmod{17},\\[3pt]
\spt(4913n+3456)&\equiv&0\pmod{17},\\[3pt]
\spt(4913n+4323)&\equiv&0\pmod{17},\\[3pt]
\spt(4913n+4612)&\equiv&0\pmod{17},\\[3pt] \label{Ono-pf-cong-2}
\spt(4913n+4901)&\equiv&0\pmod{17},\\[3pt] \label{garvan-pf-cong}
\spt(11875n+99)&\equiv&0\pmod{19},\\[3pt] \label{Ono-pf-cong-3}
\spt(12167n+9500)&\equiv&0\pmod{23},\\[3pt] \label{Ono-pf-cong-4}
\spt(24389n+806)&\equiv&0\pmod{29}.
\end{eqnarray}

All the above conjectures have been confirmed.
The congruence \eqref{garvan-pf-cong} has been proved by
Garvan \cite{Garvan-2010},
and the rest are  consequences of Theorem \ref{Ono-con}.
Indeed, Theorem \ref{Ono-con} (i) implies that if $\left(\frac{-\delta}{\ell}\right)=1$, then
\begin{equation}\label{ono-conject-prove}
\spt\left(
\frac{\ell^2(\ell n+\delta)+1}{24}
\right)
\equiv 0
\pmod{\ell}.
\end{equation}
When $\ell=11,17,23,29$,  \eqref{ono-conject-prove} becomes \eqref{Ono-pf-cong-1}$-$\eqref{Ono-pf-cong-2}, \eqref{Ono-pf-cong-3} and \eqref{Ono-pf-cong-4},  respectively.

\section{Generalizations and variations}

 In this section, we discuss  three generalizations and one variation  of the spt-function based on the relation \eqref{spt-def-eq} and three variations based on the combinatorial definition.

\subsection{The higher order spt-function of Garvan}

The first generalization of the spt-function was due to Garvan  \cite{Garvan-2011}. He defined a higher order spt-function in terms of the $k$-th symmetrized rank function and the $k$-th symmetrized crank function.

The $k$-th symmetrized rank function $\eta_k(n)$ was introduced by  Andrews  \cite{Andrews-2007}, and it is defined by
\begin{equation}\label{equ-def-etak}
\eta_{k}(n) =\sum_{m=-n}^n{m+\lfloor \frac{k-1}{2}\rfloor\choose k}N(m,n).
\end{equation}
By using $q$-identities,  Andrews  \cite{Andrews-2007} found a combinatorial interpretation of $\eta_k(n)$ in terms of $k$-marked Durfee symbol. Ji \cite{Ji-2011}  and Kursungoz \cite{Kursungoz-2011} found
combinatorial derivations of this combinatorial interpretation of $\eta_k(n)$
directly from the definition (\ref{equ-def-etak}).   When $k=2$, it is easy to check that
\[\eta_2(n)=\frac{1}{2}N_2(n),\]
where the second rank moment $N_2(n)$ is defined as in (\ref{equ-def-nkmn}).

Garvan \cite{Garvan-2011} introduced the $k$-th symmetrized crank function $\mu_k(n)$
as follows:
\begin{equation}\label{equ-def-muk}
\mu_{k}(n) =\sum_{m=-n}^n{m+\lfloor \frac{k-1}{2}\rfloor\choose k}M(m,n).
\end{equation}
A combinatorial interpretation of $\mu_k(n)$ was given by Chen, Ji and Shen \cite{Chen-Ji-Shen}. When $k=2$, it is not difficult  to derive that
\[\mu_2(n)=\frac{1}{2}M_2(n).\]

Garvan \cite{Garvan-2011}  introduced the higher order $\spt$-function $\spt_k(n)$.

\begin{defi}
For  $k\geq 1$, define
\begin{equation}\label{def-sptk}
\spt_k(n)=\mu_{2k}(n)-\eta_{2k}(n).
\end{equation}
\end{defi}

In view of \eqref{spt-def-eq}, it is easy to see  that  $\spt_k(n)$
 reduces  to  $\spt(n)$ when $k=1$.
Making use of Bailey pairs \cite{Andrews-1986}, Garvan obtained the  generating function of $\spt_k(n)$.

\begin{thm}[Garvan]\label{equ-gef-sptk-thm}
For  $k\geq 1$,
\begin{align}\label{equ-gef-sptk}
&\!\!\!\!\!\!\!\!\!\sum_{n=1}^\infty \spt_k(n)q^n\nonumber\\[3pt]
&=
\sum_{n_k\geq n_{k-1}\geq \cdots\geq n_1\geq 1}\frac{q^{n_1+n_2+\cdots +n_k}}{(1-q^{n_k})^2 (1-q^{n_{k-1}})^2\cdots (1-q^{n_{1}})^2 (q^{n_1+1};q)_\infty}.
\end{align}
\end{thm}

Setting $k=1$ in \eqref{equ-gef-sptk}, we get the generating function
\eqref{gf-spt} of $\spt(n)$. Furthermore, it can be seen from \eqref{equ-gef-sptk} that $\spt_k(n)\geq 0$ for $n,k\geq 1$. Together with  \eqref{def-sptk}, we find that
\begin{equation}\label{ine-muk-ge-etak}
\mu_{2k}(n)\geq \eta_{2k}(n).
\end{equation}

The inequality \eqref{ine-muk-ge-etak} plays a key role in the proof of an inequality  between the rank moments and the crank moments, as conjectured by Garvan \cite{Garvan-2010}.

\begin{conje}[Garvan]\label{conj-Ga-mk-ge-nk}
For  $n,k\geq 1$,
\begin{equation}\label{conj-Ga-mk-ge-nk-eq}
M_{2k}(n)\geq N_{2k}(n).
\end{equation}
\end{conje}

 Bringmann and Mahlburg \cite{Bringmann-Mahlburg-2009} showed that the above conjecture is true for $k=1,2$ and sufficiently large $n$.  For each fixed $k$, Garvan's conjecture was proved for  sufficiently large $n$ by Bringmann, Mahlburg and Rhoades \cite{Bringmann-Mahlburg-Rhoades-2011}.  Garvan \cite{Garvan-2011} confirmed  his conjecture  for all  $k$ and $n$.
He introduced  an analogue of the Stirling numbers of the second kind, denoted by $S^*(k,j)$. It is defined recursively as follows:
\begin{itemize}
\item[\rm(1)] $S^*(1,1)=1${\rm ;}
\item[\rm(2)] $S^*(k,j)=0$ if $j\leq 0$ or $j>k${\rm ;}
\item[\rm(3)] $S^*(k+1,j)=S^*(k,j-1)+j^2S^*(k,j)$ for $1\leq j\leq k+1$.
\end{itemize}
It is clear from the above recurrence relation  that  $S^*(k,j)\geq 0$.
Garvan  established the following   relations between the ordinary moments and symmetrized moments in terms of $S^*(k,j)$:
\begin{equation}\label{equ-rel-mk-muk}
M_{2k}(n)=\sum_{j=1}^k (2j)! S^*(k,j)\mu_{2j}(n)
\end{equation}
and
\begin{equation}\label{equ-rel-nk-etak}
N_{2k}(n)=\sum_{j=1}^k (2j)! S^*(k,j)\eta_{2j}(n).
\end{equation}
It follows from  \eqref{equ-rel-mk-muk} and \eqref{equ-rel-nk-etak}  that
\begin{equation}\label{equ-rel-nk-etak-rrr}
M_{2k}(n)-N_{2k}(n)=\sum_{j=1}^k (2j)! S^*(k,j)\left(\mu_{2j}(n)-\eta_{2j}(n)\right).
\end{equation}
 Invoking \eqref{ine-muk-ge-etak} we deduce that $M_{2k}(n)- N_{2k}(n)\geq 0$ for $n,k\geq1$,  and hence  Conjecture \ref{conj-Ga-mk-ge-nk} is proved.

Garvan \cite{Garvan-2011} gave a combinatorial explanation of the right-hand side of
(\ref{equ-gef-sptk}). Thus Theorem \ref{equ-gef-sptk-thm} leads to a combinatorial interpretation of $\spt_k(n)$.

\begin{thm}[Garvan]
Let $\lambda$  be a partition with $m$ different parts
\[n_1<n_2<\cdots<n_m.\]
Let $k\geq 1$,  define the weight $\omega_k(\lambda)$  of $\lambda$ as follows:
\begin{eqnarray*}
\omega_k(\lambda)=&&\sum_{m_1+\cdots+m_r=k\atop 1\leq r\leq k}{f_1+m_1-1\choose 2m_1-1}\\[3pt]
&& \qquad  \times \sum_{2\leq j_2<j_3<\cdots <j_r}{f_{j_2}+m_2\choose 2m_2}{f_{j_3}+m_3\choose 2m_3}  \cdots{f_{j_r}+m_r\choose 2m_r},
\end{eqnarray*}
where $f_j=f_j(\lambda)$ denotes the multiplicity of the part $n_j$ in $\lambda$.
Then
\[\spt_k(n)=\sum_{\lambda \in P(n)}\omega_k(\lambda).\]

\end{thm}

 Garvan \cite{Garvan-2011} also obtained     congruences of $\spt_2(n)$, $\spt_3(n)$ and $\spt_4(n)$.

\begin{thm}[Garvan]
For $n\geq 1$,
\begin{eqnarray*}
\spt_2(n)\equiv 0\pmod{5},&&\text{if } n\equiv 0,1,4\pmod{5},\\[3pt]
\spt_2(n)\equiv 0\pmod{7},&&\text{if } n\equiv 0,1,5\pmod{7},\\[3pt]
\spt_2(n)\equiv 0\pmod{11},&&\text{if } n\equiv 0\pmod{11},\\[3pt]
\spt_3(n)\equiv 0\pmod{7},&&\text{if } n\not\equiv 3,6\pmod{7},\\[3pt]
\spt_3(n)\equiv 0\pmod{2},&&\text{if } n\equiv 1\pmod{4},\\[3pt]
\spt_4(n)\equiv 0\pmod{3}, &&\text{if } n\equiv 0\pmod{3}.
\end{eqnarray*}
\end{thm}

\subsection{Generalized higher order spt-functions of Dixit and Yee}

Other generalizations of the spt-function  have been given  by Dixit and Yee   \cite{Dixit-Yee-2013}, which are based on the $j$-rank introduced by Garvan \cite{Garvan-1994}. The $j$-rank   is a generalization of Dyson's rank. For a partition $\lambda$ and $j\geq 2$, let $n_j(\lambda)$ denote the size of the $j$-th
successive Durfee square of $\lambda$,  let $c_{j}(\lambda)$ denote
the number of columns in the Ferrers diagram of  $\lambda$ with length not exceeding $ n_{j}(\lambda)$ and let $r_j(\lambda)$ denote the number of parts of $\lambda$ that lie below the $j$-th Durfee square. Then the $j$-rank of  $\lambda$ is defined  to be $c_{j-1}(\lambda)-r_{j-1}(\lambda)$. It should be noted that the $2$-rank coincides with Dyson's rank.

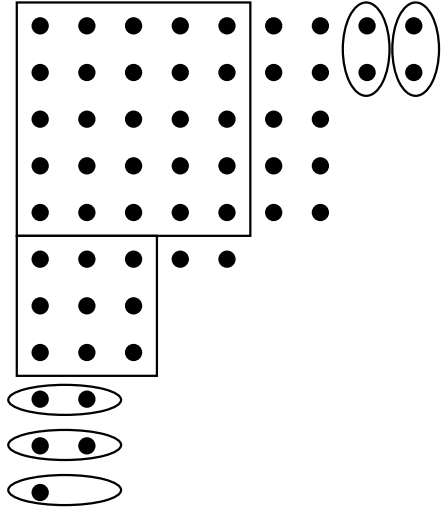
For example,   the $3$-rank of $\lambda=(9,9,7,7,7,5,3,3,2,2,1)$ is equal to $-1$, since $n_2(\lambda)=3$, $c_2(\lambda)=2$ and $r_2(\lambda)=3$, see  Figure \ref{fig5}.
\input{fig5.TpX}

Let $N_j(m,n)$ denote the number of partitions of $n$ with $j$-rank $m$.  Garvan \cite{Garvan-1994} showed that for $j\geq 2$,
\begin{equation}
\sum_{n=0}^\infty N_j(m,n)q^n=\frac{1}{(q;q)_\infty}\sum_{n=1}^\infty(-1)^{n-1}q^{\frac{n((2j-1)n-1)}{2}+|m|n}(1-q^n).
\end{equation}

Dixit and Yee \cite{Dixit-Yee-2013} defined the $j$-rank moment $_jN_k(n)$ by
\begin{equation}\label{jrankmom}
_jN_k(n) =\sum_{m=-\infty}^\infty m^k N_j(m,n).
\end{equation}
In the notation  $_jN_k(n)$, they defined  $\Spt_j(n)$ as follows.

\begin{defi}
For  $n,j \geq 1$,
\begin{equation}\label{def-sptj}
\Spt_j(n)=np(n)-\frac{1}{2}\ {_{j+1}N_2(n)}.
\end{equation}
\end{defi}
In light of \eqref{spt-moments}, it is easy to see  that  $\Spt_j(n)$
 reduces  to  $\spt(n)$ when $j=1$.

  Dixit and Yee \cite{Dixit-Yee-2013} derived the  generating function of $\Spt_j(n)$.

\begin{thm}[Dixit and Yee]
For  $j\geq 1$,
\begin{align}\label{def-sptj-gf}
&\!\!\!\!\!\!\!\!\sum_{n=1}^\infty \Spt_j(n)q^n \nonumber \\[5pt]
&=\sum_{n_j\geq 1}\sum_{n_{j-1}\geq\cdots \geq n_1\geq 0}\frac{q^{n_j}}{(1-q^{n_j})(q^{n_j};q)_\infty} {n_j\brack n_{j-1}}\cdots{n_2\brack n_1}q^{n_1^2+\cdots+n_{j-1}^2},
\end{align}
where the $q$-binomial coefficients or the Gaussian coefficients are defined by
\begin{equation}
{n\brack k}=\frac{(q;q)_n}{(q;q)_k(q;q)_{n-k}}.
\end{equation}
\end{thm}

Dixit and Yee also found a combinatorial interpretation of $\Spt_j(n)$.
To give a combinatorial explanation of the right-hand side of
\eqref{def-sptj-gf}, they introduced the $k$-th lower-Durfee square of a partition $\lambda$.
For a partition $\lambda$, take  the largest square that fits inside the Ferrers diagram of $\lambda$ starting from the lower left corner. This square is called the lower-Durfee square. If there are remaining parts above the lower-Durfee square,  then take the second lower-Durfee square in the diagram above the lower-Durfee square. Repeating this process, we are led to the third lower-Durfee square, if it exists, and so on.

The combinatorial explanation of the right-hand side of (\ref{def-sptj-gf}) also requires
a labeling of a partition, as given by Dixit and Yee.
For a partition $\lambda$, let $f_i$ denote the multiplicity of $i$ in $\lambda$.
For the $f_i$ occurrences of $i$, we label these $f_i$ parts from left to right
by $1,2,\ldots,f_i$. The labels are represented by subscripts.
For instance, $(9,8,8,8,8,6,6, 5,4,4,3)$ can be labeled as
$(9_1,8_1,8_2,8_3,8_4,6_1,6_2,5_1,4_1,4_2,3_1)$.

Using the lower-Durfee squares and the above labeling of a partition,
for a partition $\lambda$ and $j\geq 1$,  Dixit and Yee defined the weight of $\lambda$ , denoted $W_j(\lambda)$. There are two cases:

\noindent Case 1: $\lambda$  does not contain the $(j-1)$-th  lower-Durfee square. Then
$W_j(\lambda)$ is defined to be the sum of the labels of $\lambda$.

\noindent Case 2:   $\lambda$ contains the $(j-1)$-th   lower-Durfee square. Then
$W_j(\lambda)$ is defined to be the sum of labels of all the parts that are contained in and below the $(j-1)$-th   lower-Durfee square and the label of the part just  right above the $(j-1)$-th lower-Durfee square.

For example, for $\lambda=(9,8,8,8,8,6,6,5,4,4,3)$ and $j=3$, we have $W_3(\lambda)=2+3+4+1+2+1+1+2+1=17$, see Figure \ref{fig1}.
\input{fig1.TpX}

    We are now ready to state the combinatorial interpretation of $\Spt_j(n)$.

\begin{thm}[Dixit and Yee] For $j\geq 1$,
\[\Spt_j(n)=\sum_{\lambda\in P(n)}W_j(\lambda).\]
\end{thm}

Analogous to the $k$-th symmetrized rank moments $\eta_k(n)$ and the $k$-th symmetrized crank moments $\mu_k(n)$,
Dixit and Yee \cite{Dixit-Yee-2013}  defined the  $k$-th symmetrized $j$-rank function $_j\mu_{k}(n)$ by
\[_j\mu_k(n)=\sum_{m=-\infty}^\infty {m+\left\lfloor\frac{k-1}{2}\right\rfloor\choose k}N_j(m,n).\]
It can be checked that $_1\mu_k(n)=\mu_k(n)$ and $_2\mu_k(n)=\eta_k(n)$. By the
 definition \eqref{def-sptk} of the higher order spt-function $\spt_k(n)$, we see that
  \begin{equation}\label{equ-sptk-1muk-2muk}
\spt_k(n)={_1}\mu_{2k}(n)-{_2}\mu_{2k}(n).
\end{equation}
The  generalized higher order spt-function $_j\spt_k$  is defined as follows.

\begin{defi}
For $j,k\geq 1$,
\[{_j\spt}_k(n)={_j}{\mu_{2k}(n)}-{_{j+1}}\mu_{2k}(n).\]
\end{defi}

Dixit and Yee \cite{Dixit-Yee-2013} derived the  generating function of $_j\spt_k(n)$:

\begin{thm}[Dixit and Yee]
For $j,k\geq 1$,
\begin{align}\label{equ-gef-jsptk}
\sum_{n=1}^\infty {_j}\spt_k(n)q^n&=\sum_{n_k\geq \cdots\geq n_1\geq \atop m_1\geq \cdots\geq m_{j-1}\geq 1}\left(\frac{q^{n_k+\cdots+n_1}(q;q)_{n_1}}{(1-q^{n_k})^2\cdots(1-q^{n_1})^2(q^{n_1+1};q)_\infty}\right.
\nonumber\\[5pt]
& \hskip 1.27cm \times\left.\frac{q^{m_1^2+\cdots+m_{j-1}^2}} {(q;q)_{n_1-m_1}(q;q)_{m_1-m_2}
\cdots (q;q)_{m_{j-1}}}\right).
\end{align}
\end{thm}

They also gave a combinatorial explanation of the right-hand side of \eqref{equ-gef-jsptk}. Let $\lambda$ be a partition, and let $f_t$ denote the number of occurrences of $t$ in $\lambda$.  We shall use the same labeling of $\lambda$ as given before.
 For a  positive integer $k$ and a part $t$ in $\lambda$ with  label $a$,   define

\begin{align*}
g_k(\lambda,t_a)
&= {a+k-1\choose 2k-1}\\[3pt]&\qquad +
\sum_{r=2}^k \sum_{\begin{subarray}{c}
    m_1, m_2,\ldots,m_r\geq 1\\
    m_1+\cdots+m_r=k \\
    t<t_2<\cdots <t_r\leq \lambda_1
  \end{subarray}}{a+m_1-1\choose 2m_1-1}{f_{t_2}+m_2\choose 2m_2}\cdots{f_{t_r}+m_r\choose 2m_r}.
\end{align*}

\begin{defi}
For $j,k\geq 1$, define
\begin{equation}\label{equ-def-jsptk}
_j\omega_k(\lambda)=\sum_{t_a}g_k(\lambda,t_a),
\end{equation}
where the sum ranges over   the parts that are  contained in the $(j - 1)$-th lower-Durfee square except for  the last part, but also contains the part immediately above the $(j - 1)$-th lower-Durfee square.
\end{defi}

For example, let $\lambda=(5,5,5,3,3,2,2,2)$, $j=3$ and $k=2$. Label $\lambda$ as $(5_1,5_2,5_3,3_1,3_2,2_1,2_2,2_3)$. Then
\[g_2(\lambda,3_1)=0+1\cdot {3+1\choose 2}=6\]
and
\[g_2(\lambda,3_2)=1+2\cdot {3+1\choose 2}=13.\]
Moreover, from \eqref{equ-def-jsptk} we find that
\[_3\omega_2(\lambda)=g_2(\lambda,3_1)+g_2(\lambda,3_2)=6+13=19.\]
 Figure \ref{figjsptk} gives an illustration of this example.
\input{figjsptk.TpX}

  Dixit and Yee \cite{Dixit-Yee-2013} proved that $_j\spt_k(n)$ can be expressed in terms of $_j\omega_k(\lambda)$.
\begin{thm}[Dixit and Yee]
We have
\[_j\spt_k(n)=\sum_{\lambda\in P(n)}{_j}\omega_k(\lambda).\]
\end{thm}

\subsection{The ospt-function of Andrews, Chan and Kim}

A variation  of the spt-function based on relation \eqref{spt-def-eq} was given by
Andrews, Chan and Kim \cite{Andrews-Chan-Kim-2013}. In view of the symmetry properties $N(-m,n)=N(m,n)$ and $M(-m,n)=M(m,n)$, it is known that
\[N_{2k+1}(n)=M_{2k+1}(n)=0.\]

To avoid the trivial odd moments, Andrews, Chan and Kim \cite{Andrews-Chan-Kim-2013} introduced the   modified rank and crank moments $N^+_j(n)$ and $M^+_j(n)$
by considering the  unilateral sums:
\[N^+_j(n)=\sum_{m\geq 0}m^jN(m,n)\]
and
\[M^+_j(n)=\sum_{m\geq 0}m^jM(m,n).\]

They proved  the following inequality.

\begin{thm}[Andrews, Chan and Kim]\label{thm-mk+n>nk+n}
For  $n,k\geq1$,
\begin{equation}\label{thm-mk+n>nk+neq}
M^+_k(n)>N^+_k(n).
\end{equation}
\end{thm}

Bringmann and Mahlburg \cite{Bringmann-Mahlburg-2014} proved that the above inequality \eqref{thm-mk+n>nk+neq}  holds for any
fixed positive integer $k$ and sufficiently large $n$ by deriving an asymptotic
formula for $M^{+}_k(n)-N^{+}_k(n)$ stated in Theorem \ref{thm-asy-p}. When $k$ is even, this inequality \eqref{thm-mk+n>nk+neq} is equivalent to the inequality \eqref{conj-Ga-mk-ge-nk-eq} of Garvan between the rank moments and the crank moments. Chen, Ji and Zang \cite{Chen-Ji-Zang-2015} showed that the Andrews-Dyson-Rhoades conjecture \eqref{ine-conj-nsum} implies the inequality \eqref{thm-mk+n>nk+neq}.

Andrews, Chan and Kim \cite{Andrews-Chan-Kim-2013} defined the ospt-function $\ospt(n)$
as given below:

\begin{defi}\label{defi-ospt}
For $n\geq 1$,
\begin{equation}\label{defi-ospt-eq}
\ospt(n)={M}^+_1(n)-{N}^+_1(n).
\end{equation}
\end{defi}

They obtained  the   generating function of $\ospt(n)$.

\begin{thm}[Andrews, Chan and Kim]\label{thm-gef-ospt}
We have
{\small
\begin{align}\label{thm-gef-ospt-eq}
\sum_{n=0}^\infty\ospt(n)q^n&=\frac{1}{(q;q)_\infty}\sum_{i=0}^\infty \left( \sum_{j=0}^\infty q^{6i^2+8ij+2j^2+7i+5j+2}
(1-q^{4i+2})(1 -q^{4i+2 j+3})\right. \nonumber \\[5pt]
&\qquad +\left.\sum_{j=0}^\infty q^{6i^2+8ij+2j^2+5i+3j+1}(1-q^{2i+1})(1- q^{4i+2 j+2})\right).
\end{align}}
\end{thm}

Andrews, Chan and Kim  found  a combinatorial interpretation of the right-hand side of  \eqref{thm-gef-ospt-eq}, which leads to a combinatorial interpretation of $\ospt(n)$.
In doing so,  they defined even strings and odd strings of
  a partition.

  \begin{defi}
   Let $\lambda$ be a partition. A maximal consecutive sequence  $(r,r-1,\ldots,s)$ in $\lambda$ is called an even string of $\lambda$ if it satisfies the following restrictions:
   \begin{itemize}
   \item[\rm{(1)}] $r\geq 2s-2${\rm ;}
   \item[\rm{(2)}] $r$ and $s$ are even.
   \end{itemize}
Similarly, a consecutive sequence  $(r,r-1,\ldots,s)$ in $\lambda$, not necessarily maximal, is called an odd string of $\lambda$ if it satisfies the following restrictions:
\begin{itemize}
   \item[\rm{(1)}] $r+1$ is not a part of $\lambda${\rm ;}
    \item[\rm{(2)}] $s$ is odd and it appears only once in $\lambda${\rm ;}
   \item[\rm{(3)}] $r\geq 2s-1$.
   \end{itemize}

\end{defi}

For example, the partition  $\lambda=(5,4,4,3,2,2)$  contains only one odd string  $(5,4,3)$, and it does not contain any even string.
For $\lambda=(6,4,4,3,2)$, it  contains an even string $(4,3,2)$, but it does not contain any odd string.

Andrews, Chan and Kim \cite{Andrews-Chan-Kim-2013} found a combinatorial interpretation of $\ospt(n)$.

\begin{thm}[Andrews, Chan and Kim]\label{thm-com-ospt} For a partition $\lambda$, let \break ST$(\lambda)$ denote the total number of even strings and odd strings in $\lambda$.
For $n\geq 1$,
\[\ospt(n)=\sum_{\lambda \in  P(n)} \text{ST}(\lambda).\]
\end{thm}

 In light of Theorem \ref{thm-com-ospt}, Bringmann and Mahlburg \cite{Bringmann-Mahlburg-2014} proved    a monotone
 property of $\ospt(n)$ by a combinatorial argument.

\begin{thm}[Bringmann and Mahlburg]
For $n\geq 1$,
\[\ospt(n+1)\geq \ospt(n).\]
\end{thm}

They also noticed that $\ospt(n)$ and $\spt(n)$ have the same parity.
This fact can be justified as follows: Since
\[M_1^+(n)=\sum_{m\geq 0}mM(m,n)\equiv\sum_{m\geq 0}m^2M(m,n)=M_2^+(n)\pmod{2}\]
and
\[N_1^+(n)=\sum_{m\geq 0}mN(m,n)\equiv\sum_{m\geq 0}m^2N(m,n)=N_2^+(n)\pmod{2},\]
we see  that
\[\ospt(n)=M_1^+(n)-N_1^+(n)\equiv M_2^+(n)-N_2^+(n)=\spt(n)\pmod 2.\]

With the aid of the characterization of the parity  of $\spt(n)$,
Bringmann and Mahlburg  \cite{Bringmann-Mahlburg-2014}  determined the parity of $\ospt(n)$.

\begin{thm}[Bringmann and Mahlburg]
The ospt-function   ospt$(n)$ is odd if and only if $24n-1 = p^{4a+1}m^2$ for some prime $p \equiv 23  \pmod{24}$
and some integers $a, m$, where $(p,m) = 1$.
\end{thm}

Chan and Mao \cite{Chan-Mao-2014} established an upper bound and a lower bound for $\ospt(n)$,  leading to an asymptotic estimate of $\ospt(n)$.

\begin{thm}[Chan and Mao]
We have
\begin{align}
&\ospt(n)>\frac{p(n)}{4}+\frac{N(0,n)}{2}-\frac{M(0,n)}{4}\quad \text{for } n\geq 8,\label{ine-ospt-lower}\\[15pt]
&\ospt(n)<\frac{p(n)}{4}+\frac{N(0,n)}{2}-\frac{M(0,n)}{4}+\frac{N(1,n)}{2}\quad \text{for } n\geq 7,\label{ine-ospt-upper}\\[15pt]
&\ospt(n)<\frac{p(n)}{2} \quad \text{for } n\geq 3.
\end{align}
\end{thm}

An asymptotic estimate of $\ospt(n)$
can be deduced from the bounds \eqref{ine-ospt-lower} and \eqref{ine-ospt-upper}, along with
  an asymptotic property of $M(m,n)$ and $N(m,n)$ due to
 Mao \cite{Mao-2014}.

\begin{thm}[Mao] For any integer $m$, as $n\rightarrow\infty$
\begin{equation}\label{equ-asy-mln-nlm}
M(m,n)\sim N(m,n)\sim\frac{\pi}{4\sqrt{6n}}p(n).
\end{equation}
\end{thm}

By \eqref{equ-asy-mln-nlm},  we see that  as $n\rightarrow\infty$,
\[\frac{p(n)}{4}+\frac{N(0,n)}{2}-\frac{M(0,n)}{4}
\sim\frac{p(n)}{4}+\frac{N(0,n)}{2}-\frac{M(0,n)}{4}
+\frac{N(1,n)}{2}\sim\frac{1}{4}p(n).\]
Combining \eqref{ine-ospt-lower} and \eqref{ine-ospt-upper}, we arrive at
   the   asymptotic estimate \eqref{BriMaheva-for} due to Bringmann and Mahlburg
   \cite{Bringmann-Mahlburg-2014} as given in Section 5.

\subsection{The first variation of Ahlgren, Bringmann and Lovejoy}

We now turn to   three   variations of the spt-function  based on the combinatorial definition. The first variation of the spt-function was given by Ahlgren, Bringmann and Lovejoy \cite{Ahlgren-Bringmann-Lovejoy-2011}. They defined the $\M$-function  as follows.

\begin{defi}
The function  $\M(n)$  is defined to be  the total number of  smallest parts in all  partitions of $n$  without repeated odd parts and the smallest part is even.
\end{defi}

  For example,   there are  two partitions of $7$ without repeated odd parts and the smallest part is even, namely,
  \[(5,{\bf 2}),(3,{\bf 2},{\bf 2}).\]
 So we have  $\M(7)=3$.

By \cite[Section 7]{Bringmann-Lovejoy-Osburn-2010}, Ahlgren, Bringmann and Lovejoy \cite{Ahlgren-Bringmann-Lovejoy-2011} derived the  generating function of $\M(n)$.

\begin{thm}[Ahlgren, Bringmann and Lovejoy]
We have
{\small
\begin{eqnarray}\label{equ-gf-m2spt}
&&\sum_{n=1}^\infty\M(n)q^n=\frac{(-q;q^2)_\infty}{(q^2;q^2)_\infty}\nonumber \\[3pt]
&&\hskip 3cm \times \left(\sum_{n=1}^\infty\frac{nq^{2n}}{1-q^{2n}}+\sum_{n=-\infty\atop n\neq 0}^{\infty} \frac{(-1)^nq^{2n^2+n}}{(1-q^{2n})^2}\right).
\end{eqnarray}}
\end{thm}

Jennings-Shaffer \cite{Jennings-2015} showed that the function $\M(n)$ can be expressed as the difference between the
 symmetrized $M_2$-rank moments   and the symmetrized residue crank moments of  partitions without repeated odd parts.  Let us first recall the definitions of the $M_2$-rank of a partition without repeated odd parts   and the residue crank of a partition without repeated odd parts.

  Let $\lambda$ be a partition  without repeated odd parts, the $M_2$-rank of $\lambda$ was defined by Berkovich and Garvan \cite{Berkovich-Garvan-2002} as stated below:
\begin{equation}\label{equ-def-m2-rank-ord}
M_2\text{-rank}(\lambda)=\left\lceil\frac{\lambda_1}{2}
\right\rceil-l(\lambda).
\end{equation}

 The residue crank of $\lambda$ was defined by  Garvan and Jennings-Shaffer \cite{Garvan-Jennings-2014} which is related to the crank of an ordinary partition. Let $\lambda=(\lambda_1,\lambda_2,\ldots,\lambda_l)$ be a partition without repeated odd parts, define  $\lambda^e$ to be the ordinary partition obtained from  $\lambda$ by omitting  odd parts of $\lambda$ and dividing each even part by $2$.  The residue crank of $\lambda$  is defined to be  the crank of $\lambda^e$.

 For example, let $\lambda=(11,7,6,5,4,4,3,2,2)$, then $\lambda_1=11$, $l(\lambda)=9$ and $\lambda^e=(3,2,2,1,1)$. Hence the $M_2$-rank of $\lambda$ is equal to $-3$ and the residue crank of $\lambda$ is equal to the crank of $\lambda^e$, which equals $-1$.

Let $N2(m,n)$ denote the number of partitions of $n$ without repeated odd parts  such that  $M_2$-rank is equal to $m$. Let $M2(m,n)$ denote the number of partitions of $n$ without repeated odd parts such that  the residue crank is equal to $m$. The $k$-th symmetrized $M_2$-rank moments $\eta 2_k(n)$ and the $k$-th symmetrized residue crank moments $\mu2_k(n)$  of partitions without repeated odd parts were defined by  Jennings-Shaffer \cite{Jennings-2015}  as follows:
\begin{align*}
\eta 2_k(n)&=\sum_{m=-\infty}^\infty {m+\left\lfloor\frac{k-1}{2}\right\rfloor\choose k} N2(m,n),\\[5pt]
\mu 2_k(n)&=\sum_{m=-\infty}^\infty {m+\left\lfloor\frac{k-1}{2}\right\rfloor\choose k} M2(m,n).
\end{align*}

Analogue to the relation \eqref{spt-def-eq} for $\spt(n)$,
Jennings-Shaffer \cite{Jennings-2015} established the following connection.

\begin{thm}[Jennings-Shaffer]
For   $n\geq 1$,
\begin{equation}\label{m2spt-function}
\M(n)=\mu2_2(n)-\eta2_2(n).
\end{equation}
\end{thm}

 The following congruences of $\M(n)$ mod $3$ and $5$ were given by Garvan and Jennings-Shaffer \cite{Garvan-Jennings-2014}.
 \begin{thm}[Garvan and Jennings-Shaffer]
 For $n\geq0$,
 \begin{eqnarray*}
\M(3n+1)&\equiv& 0 \pmod 3,\\[3pt]
\M(5n+1) &\equiv& 0 \pmod 5,\\[3pt]
\M(5n+3) &\equiv& 0 \pmod {5}.
\end{eqnarray*}
\end{thm}

 Jennings-Shaffer \cite{Jennings-2015-1} provided alternative proofs of the above congruences. Furthermore, he showed that
\begin{thm}[Jennings-Shaffer]
For $n\geq0$,
\begin{eqnarray*}
\M(27n+26) &\equiv& 0 \pmod 5, \\[3pt]
\M(125n+97) &\equiv& 0 \pmod 5, \\[3pt]
\M(125n+122) &\equiv& 0 \pmod 5.
\end{eqnarray*}
\end{thm}

 Ahlgren, Bringmann and Lovejoy \cite{Ahlgren-Bringmann-Lovejoy-2011} established    Ramanujan-type congruences of $\M(n)$ modulo   powers of $\ell$ for any prime $\ell\geq 3$.

\begin{thm}[Ahlgren, Bringmann and Lovejoy]
Let  $\ell\geq 3$  be a prime, and let $ m,n\geq 1$.
\begin{itemize}
\item[\rm{(i)}]  If $\left(\frac{-n}{\ell}\right)=1$, then
\[\M\left((\ell^{2m}n+1)/8\right)\equiv 0\pmod{\ell^m}.\]
\item[\rm{(ii)}]
\[\M\left((\ell^{2m+1}n+1)/8\right)\equiv \left(\frac{2}{\ell}\right)\M\left((\ell^{2m-1}n+1)/8\right)
\pmod{\ell^m}.\]
\end{itemize}
\end{thm}

 Hecke-type congruences   of $\M(n)$ mod $2,2^2, 2^3,3$ and $5$ have been found by
Andersen   \cite{Andersen-2013}.

\begin{thm}[Andersen]\label{thm-heck-m2}
 Let $\ell\geq3$ be a prime. Define $s_\ell=(\ell^2-1)/8$ and
 \[ \beta=\begin{cases}
 1,&\text{if }\ell\equiv 3\pmod{8},\\[3pt]
 2,&\text{if }\ell\equiv 5\pmod{8},\\[3pt]
 3,&\text{if }\ell\equiv 1,7\pmod{8}.
 \end{cases}\]
 For $t\in\{2^\beta,3,5\},\ell\neq t$ and $n\geq 1$,
 \begin{gather*}
 \M(\ell^2n-s_\ell)+\left(\frac{2}{\ell}\right)\left(\frac{1-8n}{\ell}\right)
 \M(n)+\ell \M\left((n+s_\ell)/\ell^2\right) \\[3pt]
 \equiv\left(\frac{2}{\ell}\right)(1+\ell)
 \M(n)\pmod{t}.
 \end{gather*}
 \end{thm}

In analogy with  the higher order $\spt$-function $\spt_k(n)$, Jennings-Shaffer \cite{Jennings-2015-1} defined the  higher order function $\M_k(n)$ in terms of  the $k$-th symmetrized $M_2$-rank moments $\eta2_k(n)$   and the $k$-th symmetrized residue crank moments $\mu2_k(n)$  for partitions without repeated odd parts.

\begin{defi}For $k\geq 1$, define
\begin{equation*}
\M_k(n)=\mu2_{2k}(n)-\eta2_{2k}(n).
\end{equation*}
\end{defi}

Using \eqref{m2spt-function}, it is clear to see that $\M_k(n)$ reduces to $\M(n)$ when $k=1$.
Jennings-Shaffer \cite{Jennings-2015} also obtained the   generating function of $\M_k(n)$.

\begin{thm}[Jennings-Shaffer]
We have
\begin{align}\label{equ-gef-m2spt-k}
&\!\!\!\!\sum_{n=1}^\infty \M_k(n)q^n\nonumber \\[3pt]
&=\sum_{n_k\geq n_{k-1}\geq \cdots\geq n_1\geq 1} \frac{(-q^{2n_1+1};q^2)_\infty q^{2n_1+2n_2+\cdots+2n_k}}{(1-q^{2n_k})^2(1-q^{2n_{k-1}})^2\cdots (1-q^{2n_1})^2 (q^{2n_1+2};q^2)_\infty}.
\end{align}
\end{thm}

By interpreting  the right-hand side of \eqref{equ-gef-m2spt-k} combinatorially, Jennings-Shaffer \cite{Jennings-2015} found a combinatorial interpretation of $\M_k(n)$.
Let $P_o(n)$ denote the set of partitions of $n$ without repeated odd parts and the smallest part is even. For a partition $\lambda\in P_o(n)$, assume that there are $r$ different even parts in $\lambda$, namely,
\[2t_1<2t_2<\cdots<2t_r.\]
Let $f_j=f_j(\lambda)$ denote the frequency of the part $2t_j$ in $\lambda$.
For a fixed integer $k\geq 1$, Jennings-Shaffer \cite{Jennings-2015} defined
$\omega_k(\lambda)$ as follows:
\begin{equation}\label{equ-def-omega-k-l}
\omega_k(\lambda)=\sum_{m_1+m_2+\cdots+m_s=k\atop
1\leq s\leq k}{f_1+m_1-1\choose 2m_1-1}\times
\sum_{2\leq j_2<j_3<\cdots<j_s}\prod_{i=2}^s{f_{j_i}+m_i\choose 2m_i}.
\end{equation}

For example, let $k=2$ and   $\lambda=(10,10,9,5,4,3,2,2,2)$ be a partition in $P_o(47)$,  there are three distinct even parts in $\lambda$.
Thus $r=3$,  $f_1=3$, $f_2=1$ and $f_3=2$.
By the definition \eqref{equ-def-omega-k-l} of $\omega_k(\lambda)$, we have
$\omega_2(\lambda)=16.$

With the above notation, Jennings-Shaffer \cite{Jennings-2015} found a combinatorial interpretation of $\M_k(n)$.

\begin{thm}[Jennings-Shaffer]
For $n\geq 1$,
\begin{equation}\label{equ-def-eta2-2kn}
\M_k(n)=\sum_{\lambda\in P_o(n)} \omega_k(\lambda).
\end{equation}
\end{thm}

Jennings-Shaffer \cite{Jennings-2015-1} also obtained the following congruences of \break $\M_2(n)$.

\begin{thm}[Jennings-Shaffer]
For $n\geq1$,
\begin{eqnarray*}
\M_2(n)&\equiv& 0\pmod{3},~\text{if}~n\equiv0\pmod9,\\[3pt]
\M_2(n)&\equiv& 0\pmod{5},~\text{if}~n\equiv0\pmod5,\\[3pt]
\M_2(n)&\equiv& 0\pmod{5},~\text{if}~n\equiv1\pmod5,\\[3pt]
\M_2(n)&\equiv& 0\pmod{5},~\text{if}~n\equiv3\pmod5.
\end{eqnarray*}
\end{thm}

\subsection{The second variation of Bringmann, Lovejoy and Osburn}

The second variation of the spt-function was due to Bringmann, Lovejoy and Osburn \cite{Bringmann-Lovejoy-Osburn-2009}, which is defined on  overpartitions. Recall that Corteel and Lovejoy \cite{Corteel-Lovejoy-2004} defined an overpartition of $n$ as a partition of $n$ in
which the first occurrence of a part may be overlined.  Bringmann, Lovejoy and Osburn \cite{Bringmann-Lovejoy-Osburn-2009} introduced    three spt-type functions.

\begin{defi}[Bringmann, Lovejoy and Osburn]
\

\begin{itemize}
\item[\rm(1)] The function $\overline{\spt}(n)$  is defined to be the total number of smallest parts in all overpartitions  of $n$.

\item[\rm(2)] The function $\overline{\spt1}(n)$ is defined to be the total number of smallest parts in all overpartitions  of $n$ with  the smallest part being odd.

\item[\rm(3)] The function $\overline{\spt2}(n)$ is defined to be  the total number of smallest parts in all overpartitions  of $n$  with the smallest part being even.
\end{itemize}
\end{defi}

For example, there are $14$ overpartitions of $4$:
\[\begin{array}{lllllllllllll}
(4)&(\bar{4})&(3,1)&(\bar{3},1)&(3,\bar{1})&(\bar{3},\bar{1})&(2,2),\\[5pt]
(\bar{2},2)&(2,1,1)&(\bar{2},1,1)&(2,\bar{1},1)&(\bar{2},\bar{1},1)& (1,1,1,1)&(\bar{1},1,1,1).
\end{array}
\]
We have $\overline{\spt}(4)=26$, $\overline{\spt1}(4)=20$ and $\overline{\spt2}(4)=6$.

Analogous  to the relation \eqref{spt-def-eq} for the spt-function,  the functions $\overline{\spt}(n)$ and $\overline{\spt2}(n)$ can also be expressed as the differences of the rank and the crank moments of    overpartitions.  The definitions of the  rank and the crank moments of   overpartitions are based on the two definitions of the  rank of an overpartition  and
the two definitions of the crank of an overpartition.
Although there are four possibilities, only two of them have been studied.

For an overpartition $\lambda$, there are two kinds of ranks. One is called the D-rank  introduced by Lovejoy \cite{Lovejoy-2005} and the other is called the $M_2$-rank  introduced by Lovejoy \cite{Lovejoy-2008}. The  D-rank of $\lambda$ is  defined as the largest part  minus the number of  parts.  To define the $M_2$-rank,
let $\lambda_o$ denote the  partition consisting
of  non-overlined odd parts of $\lambda$. Then $M_2\text{-rank}(\lambda)$ can be defined as follows:
\[M_2\text{-rank}(\lambda)=\left\lceil\frac{\lambda_1}{2}\right\rceil-l(\lambda)
+l(\lambda_o)-\chi(\lambda),\]
where $\chi(\lambda)=1$ if the largest part of $\lambda$ is odd and non-overlined and $\chi(\lambda)=0$ otherwise.

For example, for an overpartition $\lambda=(\bar{9},9,7,\bar{6},5,5,\bar{4},3,2,\bar{1})$, we see that $\text{D-rank}(\lambda)=9-10=-1$. Moreover, since  $\lambda_o=(9,7,5,5,3)$ and $\chi(\lambda)=0$, we have
$M_2\text{-rank}(\lambda)=0.$

Bringmann, Lovejoy and Osburn \cite{Bringmann-Lovejoy-Osburn-2009} defined the first and second residue crank of an overpartition.
The first residue crank of an overpartition is
defined as the crank of the partition consisting of  non-overlined parts. The second
residue crank is defined as the crank of the subpartition consisting of all   the even non-overlined
parts divided by two.

For example, for   $\lambda=(\bar{9},9,7,\bar{6},5,5,\bar{4},4,3,2,\bar{1})$, the
 partition consisting of non-overlined parts of $\lambda$ is $(9,7,5,5,4,3,2)$.
The first residue crank of $\lambda$ is $9$. The partition formed by even non-overlined parts of $\lambda$ is  $(4,2)$. So the second residue crank of $\lambda$ is equal to the crank of $(2,1)$, which is equal to $0$.

We are now in a position to present the definitions of the rank and the crank moments of  overpartitions. Let $\overline{N}(m,n)$ denote the number of overpartitions of $n$ with the D-rank
$m$, and let $\overline{N2}(m,n)$ denote the number of overpartitions of $n$ with the $M_2$-rank   $m$.  Notice that there are two kinds of ranks of overpartitions. Consequently, there are two possibilities to define the rank moments of overpartitions. The two rank moments are defined as follows:
\begin{eqnarray}
\overline{N}_k(n)&=&\sum_{m=-\infty}^\infty m^k\overline{N}(m,n),\\[3pt]
\overline{N2}_k(n)&=&\sum_{m=-\infty}^\infty m^k\overline{N2}(m,n).
\end{eqnarray}

Similarly,  let $\overline{M}(m,n)$  denote the number of overpartitions of $n$ with the first residue crank   $m$ and let $\overline{M2}(m,n)$  denote the number of overpartitions of $n$ with the second residue crank $m$.  The two crank moments are defined by
\begin{eqnarray}
\overline{M}_k(n)&=&\sum_{m=-\infty}^\infty m^k\overline{M}(m,n),\\[3pt]
\overline{M2}_k(n)&=&\sum_{m=-\infty}^\infty m^k\overline{M2}(m,n).
 \end{eqnarray}

Bringmann, Lovejoy and Osburn   \cite{Bringmann-Lovejoy-Osburn-2009} deduced the following relations on $\overline{\spt}(n)$ and $\overline{\spt2}(n)$.

\begin{thm}[Bringmann, Lovejoy and Osburn]\label{thm-ngf-sptbar}
For $n\geq 1$,
\begin{eqnarray}
\overline{\spt}(n)&=&\overline{M}_2(n)-\overline{N}_2(n), \\[5pt]
\overline{\spt2}(n)&=&\overline{M2}_2(n)-\overline{N2}_2(n).
\end{eqnarray}
\end{thm}

In light of Theorem \ref{thm-ngf-sptbar}, Bringmann, Lovejoy and Osburn \cite{Bringmann-Lovejoy-Osburn-2009} proved the following congruences:

\begin{thm}[Bringmann, Lovejoy and Osburn]\label{equ-cong-ovep-1}
For $n\geq 1$,
\begin{eqnarray}
\overline{\spt2}(n)\equiv\overline{\spt2}(n)&\equiv& 0\pmod{3},~\text{if}~n\equiv0,1\pmod3,\label{cong-op-1}\\[5pt]
\overline{\spt}(n)&\equiv&0\pmod{3},~\text{if}~n\equiv0\pmod3,
\label{cong-op-2}\\[5pt]
\overline{\spt2}(n)&\equiv&0\pmod{5},
~\text{if}~n\equiv3\pmod5,\label{cong-op-3}\\[5pt]
\overline{\spt1}(n)&\equiv&0\pmod{5},~\text{if}~n\equiv0\pmod5.\label{cong-op-4}
\end{eqnarray}
Moreover, if $\ell\geq 5$ is a prime, then the following congruence holds:
\[\overline{\spt1}(\ell^2 n)+\left(\frac{-n}{\ell}\right)\overline{\spt1}(n) +\ell\overline{\spt1}\left(n/\ell^2\right)\equiv
(\ell+1)\overline{\spt1}(n)\pmod{3}.\]
\end{thm}

An alternative proof of  the congruence  \eqref{cong-op-1}  was given by Jennings-Shaffer  \cite{Jennings-p-4}. The combinatorial interpretations of the congruences \eqref{cong-op-1}--\eqref{cong-op-4} were given by  Garvan and Jennings-Shaffer \cite{Garvan-Jennings-2014}.
Ahlgren, Bringmann and Lovejoy \cite{Ahlgren-Bringmann-Lovejoy-2011} derived  Ramanujan-type congruences of $\overline{\spt1}(n)$ modulo   powers of a prime $\ell$, which are similar to the Ramanujan-type congruences of $\spt(n)$ modulo  powers of a prime $\ell$.

\begin{thm}[Ahlgren, Bringmann and Lovejoy]
Let  $\ell\geq 3$ be a prime, and let $m, n\geq 1$.
\begin{itemize}
\item[\rm{(1)}] If $(\frac{-n}{\ell})=1$, then
\[\overline{\spt1}(\ell^{2m} n)\equiv 0\pmod{\ell^m}.\]
\item[\rm{(2)}]
\[\overline{\spt1}(\ell^{2m+1}n)\equiv \overline{\spt1}(\ell^{2m-1}n)\pmod{\ell^m}.\]
\end{itemize}
\end{thm}

Andersen  \cite{Andersen-2013} obtained  Hecke-type congruences of  $\overline{\spt1}(n)$ mod $2^6$, $2^7$, $2^8$, $3$ and $5$.
\begin{thm}[Andersen]
 Let $\ell\geq3$ be a prime, and define
 \[\beta=\begin{cases}
 6,&\text{if }\ell\equiv 3\pmod{8},\\[3pt]
 7,&\text{if }\ell\equiv 5,7\pmod{8},\\[3pt]
 8,&\text{if }\ell\equiv 1\pmod{8}.
 \end{cases}\]
 For $t\in\{2^\beta,3,5\}$, $\ell\neq t$ and $n\geq 1$,
 \[\overline{\spt1}(\ell^2n)+\left(\frac{-n}{\ell}\right)\overline{\spt1}(n)+
 \ell\overline{\spt1}(n/\ell^2)\equiv (1+\ell)\overline{\spt1}(n)\pmod{t}.\]
 \end{thm}

It is readily seen
 that $\overline{\spt1}(n)$, $\overline{\spt2}(n)$ and $\overline{\spt}(n)$   are all  even.
 Congruences   of these  functions modulo $4$ were investigated by Garvan and Jennings-Shaffer \cite{Garvan-Jennings-2014}.

\begin{thm}[Garvan and Jennings-Shaffer]\label{thm-op-par}
For $n\geq 1$,
\begin{itemize}
\item[{\rm (1)}]$\overline{\spt}(n)\equiv 2\pmod{4}$ if and only if $n$ is a square or twice a square{\rm ;}
\item[{\rm (2)}] $\overline{\spt1}(n)\equiv 2\pmod{4}$ if and only if $n$ is an odd square{\rm ;}

\item[{\rm (3)}]$\overline{\spt2}(n)\equiv 2\pmod{4}$ if and only if $n$ is an even
square or twice a square.
\end{itemize}
\end{thm}

Moreover, they introduced a statistic $\overline{\sptcrank}$ defined on a marked overpartition, which leads to   combinatorial interpretations of  the above congruences.

The following recurrence relation  of $\overline{\spt1}(n)$ was given by Ahlgren and Andersen \cite{Ahlgren-Andersen-2015}.

\begin{thm}[Ahlgren and Andersen]
Let
$$s(n) =\sum_{d|n}\min\left(d,\frac{n}{d}\right).$$
For  $n>0$,
\[\sum_k (-1)^k\overline{\spt1}(n-k^2)=b(n),\]
where
\begin{equation*}
b(n) =\begin{cases}
2s(n), &\text{if }n\text{ is odd},\\[3pt]
-4s(n/4), &\text{if }n\equiv 0\pmod{4},\\[3pt]
0, &\text{if }n\equiv 2\pmod{4}.
\end{cases}
\end{equation*}
\end{thm}

In view of the symmetry properties $\overline{N}(-m,n)=\overline{N}(m,n)$ and $\overline{M}(-m,n)$\break $=\overline{M}(m,n)$, we see that
\[\overline{N}_{2k+1}(n)=\overline{M}_{2k+1}(n)=0.\]
Similarly, to avoid the trivial odd moments,  Andrews, Chan, Kim and Osburn  \cite{Andrews-Chan-Kim-Osburn-2013} introduced the modified rank and crank moments $\overline{N}_k^+(n)$ and $\overline{M}_k^+(n)$ for overpartitions:
\begin{equation}
\overline{N}_k^+(n)=\sum_{m\geq 1}m^k\overline{N}(m,n)
\end{equation}
and
\begin{equation}
\overline{M}_k^+(n)=\sum_{m\geq 1}m^k\overline{M}(m,n).
\end{equation}

They defined the ospt-function $\overline{\ospt}(n)$ for overpartitions which is  in the spirit of the ospt-function $\ospt(n)$ for ordinary partitions.
\begin{defi}For $n\geq 1$,
\begin{equation}\label{defi-ospt-b-eqn}
\overline{\ospt}(n)=\overline{M}_1^+(n)-\overline{N}_1^+(n).
\end{equation}
\end{defi}

Andrews, Chan, Kim and Osburn \cite{Andrews-Chan-Kim-Osburn-2013}
defined  even strings and   odd strings of an overpartition, and provided a combinatorial interpretation of $\overline{\ospt}(n)$.

Jennings-Shaffer \cite{Jennings-2015} defined the higher order $\spt$-functions for  overpartitions
by using the $k$-th symmetrized rank and crank moments for  overpartitions.
There are two symmetrized rank   moments for overpartitions:
\begin{equation}\label{equ-def-etak-op}
\overline{\eta}_{k}(n)=\sum_{m=-n}^n{m+\lfloor \frac{k-1}{2}\rfloor\choose k}\overline{N}(m,n)
\end{equation}
and
\begin{equation}\label{equ-def-etak-op2}
\overline{\eta2}_{k}(n)=\sum_{m=-n}^n{m+\lfloor \frac{k-1}{2}\rfloor\choose k}\overline{N2}(m,n).
\end{equation}
There are also two  symmetrized crank   moments for overpartitions:
\begin{equation}\label{equ-def-etak-op3}
\overline{\mu}_{k}(n)=\sum_{m=-n}^n{m+\lfloor \frac{k-1}{2}\rfloor\choose k}\overline{M}(m,n)
\end{equation}
and
\begin{equation}\label{equ-def-etak-op4}
\overline{\mu2}_{k}(n)=\sum_{m=-n}^n{m+\lfloor \frac{k-1}{2}\rfloor\choose k}\overline{M2}(m,n).
\end{equation}
The two  higher order $\spt$-functions for overpartitions are defined as follows.

\begin{defi}
For  $k\geq 1$,
\begin{equation}
\overline{\spt}_k(n)=\overline{\mu}_{2k}(n)-\overline{\eta}_{2k}(n),
\end{equation}
\begin{equation}
\overline{\spt2}_k(n)=\overline{\mu2}_{2k}(n)-\overline{\eta2}_{2k}(n).
\end{equation}
\end{defi}

Using   Bailey pairs, Jennings-Shaffer \cite{Jennings-2015} obtained the   generating functions of $\overline{\spt}_k(n)$ and $\overline{\spt2}_k(n)$.

\begin{thm}[Jennings-Shaffer]
For  $k\geq 1$,
\begin{align}\label{gf-sptk}
&\!\!\!\!\!\!\sum_{n=1}^\infty \overline{\spt}_k(n)q^n\nonumber\\[3pt]&=\sum_{n_k\geq n_{k-1}\geq \cdots\geq n_1\geq 1}\frac{q^{n_1+n_2+\cdots+n_k}(-q^{n_1+1};q)_\infty}{(1-q^{n_k})^2
(1-q^{n_{k-1}})^2\cdots(1-q^{n_1})^2(q^{n_1+1};q)_\infty},
\end{align}
\begin{align}\label{gf-barsptk}
&\!\!\!\!\sum_{n=1}^\infty \overline{\spt2}_k(n)q^n\nonumber\\[3pt]&=\sum_{n_k\geq n_{k-1}\geq \cdots\geq n_1\geq 1}\frac{q^{2n_1+2n_2+\cdots+2n_k}(-q^{2n_1+1};q)_\infty}{(1-q^{2n_k})^2
(1-q^{2n_{k-1}})^2\cdots(1-q^{2n_1})^2(q^{2n_1+1};q)_\infty}.
\end{align}
\end{thm}

By interpreting the right-hand sides  of \eqref{gf-sptk} and \eqref{gf-barsptk}
 based on vector partitions, Jennings-Shaffer found   combinatorial explanations of $\overline{\spt}_k(n)$ and $ \overline{\spt2}_k(n)$.

\subsection{The third variation of  Andrews, Dixit and Yee}

The third variation of the spt-function was introduced by Andrews, Dixit and Yee \cite{Andrews-Dixit-Yee-2015}. Let $p_\omega(n)$ denote the number of partitions of $n$ in which each odd part is less than twice the smallest part.   They defined $\spt_\omega(n)$ as follows.

\begin{defi}
 The function $\spt_\omega(n)$  is defined to be the number of smallest parts in the partitions enumerated by $p_\omega(n)$.
 \end{defi}

  For example,  for $n=4$,  there are  four partitions  counted by $p_\omega(4)$, namely,
 \[(4)\quad(2,2)\quad(2,1,1)\quad(1,1,1,1).\]
 We have $p_\omega(4)=4$ and $\spt_\omega(4)=9$.

 They derived the   generating function of $\spt_\omega(n)$.

\begin{thm}[Andrews, Dixit and Yee] We have
\begin{align}\label{equ-gf-spt-omega}
&\!\!\!\!\!\sum_{n=1}^\infty\spt_\omega(n)q^n\nonumber\\[3pt]
&=\frac{1}{(q^2;q^2)_\infty} \sum_{n=1}^\infty\frac{nq^n}{1-q^n}+ \frac{1}{(q^2;q^2)_\infty}\sum_{n=1}^\infty \frac{(-1)^n(1+q^{2n})q^{n(3n+1)}}
{(1-q^{2n})^2}.
\end{align}
\end{thm}

Using the above  generating function, Andrews, Dixit and Yee \cite{Andrews-Dixit-Yee-2015} proved the following congruences of $\spt_\omega(n)$.

 \begin{thm}[Andrews, Dixit and Yee] For $n\geq 0$,
\begin{eqnarray}\label{Andrews-Dixit-Yee-2015-1}
\spt_\omega(5n+3)&\equiv &0\pmod{5},\\[3pt]
\spt_\omega(10n+7)&\equiv &0\pmod{5},\\[3pt]
\spt_\omega(10n+9)&\equiv &0\pmod{5}.
\end{eqnarray}
\end{thm}

Employing the generating function \eqref{equ-gf-spt-omega}, Wang \cite{Wang-2017} derived the  generating function of $\spt_\omega(2n+1)$.

\begin{thm}[Wang]
We have
\begin{equation}\label{equ-gef-spt-omega-2n+1}
\sum_{n=0}^\infty \spt_\omega (2n+1)q^n=\frac{(q^2;q^2)_\infty^8}{(q;q)_\infty^5}.
\end{equation}
\end{thm}

Wang \cite{Wang-2017} also posed  two  conjectures on   congruences of $\spt_\omega(n)$ modulo
arbitrary powers of $5$.

\begin{conje}[Wang]
For   $k\geq 1$ and $n\geq 0$,
\begin{eqnarray*}
\spt_\omega\left(2\cdot 5^{2k-1}n+\frac{7 \cdot 5^{2k-1}+1}{12}\right)&\equiv& 0\pmod{5^{2k-1}}.
\end{eqnarray*}
\end{conje}

\begin{conje}[Wang]
For   $k\geq 1$ and $n\geq 0$,
\begin{eqnarray*}
\spt_\omega\left(2\cdot 5^{2k}n+\frac{11 \cdot 5^{2k}+1}{12}\right)&\equiv& 0\pmod{5^{2k}}.
\end{eqnarray*}
\end{conje}

Jang and Kim \cite{Jang-Kim-2016} obtained a congruence  of $\spt_{\omega}(n)$
via the mock modularity of its generating function.

\begin{thm}[Jang and Kim]
Let $\ell \geq5$ be a prime, and let $j,m$ and $n$ be positive integers with $\left(\frac{n}{\ell}\right)=-1$. If $m$ is sufficiently large, then there are infinitely many primes $Q\equiv-1\pmod{576\ell ^j}$ satisfying
\begin{equation}
\spt_{\omega}\left(\frac{Q^3\ell^mn+1}{12}\right)\equiv0\pmod{\ell^j}.
\end{equation}
\end{thm}

An overpartition analogue of the function $\spt_\omega(n)$ was defined by Andrews, Dixit, Schultz and Yee \cite{Andrews-Dixit-Schultz-Yee}.

\begin{defi}
The function $\overline{\spt}_\omega(n)$ is defined to be the number of smallest parts in the overpartitions of $n$ in which the smallest part is always overlined and all odd parts are less than twice the smallest part.
\end{defi}

They obtained the generating function of $\overline{\spt}_\omega(n)$.

\begin{thm}[Andrews, Dixit, Schultz and Yee]
We have
\begin{equation}\label{equ-gf-spt-ome-ove}
\sum_{n=1}^\infty \overline{\spt}_\omega(n)q^n= \frac{(-q^2;q^2)_\infty}{(q^2;q^2)_\infty}
\sum_{n=1}^\infty \frac{nq^n}{(1-q^n)}+2\frac{(-q^2;q^2)_\infty} {(q^2;q^2)_\infty}\sum_{n=1}^\infty \frac{(-1)^n q^{2n(n+1)}}{(1-q^{2n})^2}.
\end{equation}
\end{thm}

Based on the generating function \eqref{equ-gf-spt-ome-ove}, they derived the following congruences of $\overline{\spt}_\omega(n)$ mod 3, 5 and 6.

\begin{thm}[Andrews, Dixit, Schultz and Yee]
For $n\geq0$,
\begin{eqnarray*}
\overline{\spt}_\omega(3n)&\equiv &0\pmod{3},\\[3pt]
\overline{\spt}_\omega(3n+2)&\equiv &0\pmod{3},\\[3pt]
\overline{\spt}_\omega(10n+6)&\equiv &0\pmod{5},\\[3pt]
\overline{\spt}_\omega(6n+5)&\equiv &0\pmod{6}.
\end{eqnarray*}
\end{thm}

They also characterized  the parity of $\overline{\spt}_\omega(n)$.

\begin{thm}[Andrews, Dixit, Schultz and Yee]
For $n\geq1$, $\overline{\spt}_\omega(n)$ is odd if and only if $n=k^2$ or $2k^2$ for some $k\geq 1$.
\end{thm}

Moreover, they found a congruence of $\overline{\spt}_\omega(n)$ modulo  $4$.

\begin{thm}[Andrews, Dixit, Schultz and Yee]
For $n\geq1$,
\[\overline{\spt}_\omega(7n)\equiv \overline{\spt}_\omega(n/7)\pmod{4},\]
where we adopt the convention that $\overline{\spt}_\omega(x)=0$ if $x$ is not a positive integer.
\end{thm}

By \eqref{equ-gf-spt-omega}, Wang \cite{Wang-2017} obtained  the   generating function of $\overline{\spt}_\omega(2n+1)$.

 \begin{thm}[Wang]
 We have
 \begin{equation}\label{equ-gef-over-spt-omega-2n+1}
 \sum_{n=0}^\infty \overline{\spt}_\omega(2n+1)q^n=\frac{(q^2;q^2)_\infty^9}{(q;q)_\infty^6}.
 \end{equation}
 \end{thm}

In light of \eqref{equ-gef-over-spt-omega-2n+1}, Wang derived the following congruences of $\overline{\spt}_\omega(n)$.

\begin{thm}[Wang]
For   $n\geq 0$,
\begin{eqnarray*}
\overline{\spt}_\omega(8n+7)&\equiv &0\pmod{4},\\[3pt]
\overline{\spt}_\omega(6n+5)&\equiv &0\pmod{9},\\[3pt]
\overline{\spt}_\omega(18n+r)&\equiv &0\pmod{9},\qquad \text{for } r=9 \text{ or } 15,\\[3pt]
\overline{\spt}_\omega(22n+r)&\equiv &0\pmod{11},\qquad \text{for }  r=7,11,13,17,19,\text{ or }21,\\[3pt]
\overline{\spt}_\omega(162n+r)&\equiv &0\pmod{27},\qquad \text{for } r=81 \text{ or }135.
\end{eqnarray*}
\end{thm}

There are   other variations of the spt-function, and we just mention the main ideas of these variations.  Jennings-Shaffer \cite{Jennings-p-1,Jennings-p-2,Jennings-p-3} introduced  several $\spt$-type functions arising from Bailey pairs and derived several   Ramanujan-type congruences.
 Garvan and Jennings-Shaffer \cite{Garvan-Jennings-2016} discovered more $\spt$-type functions  and found some  congruences of these $\spt$-type functions.
Patkowski  \cite{Patkowski,Patkowski-1,Patkowski-2}  also defined several $\spt$-type functions based on Bailey pairs. Furthermore, Patkowski obtained   generating functions and congruences of these functions. Sarma, Reddy, Gunakala and Comissiong \cite{Sarma-Reddy-Gunakala-Comissiong-2011} defined a more general function,
 in the notation $\spt_i(n)$, as the total number of the $i$-th smallest part in all partitions of $n$.

\section{Asymptotic properties}

In this section, we present   asymptotic formulas for the spt-function and its
 variations.
By applying the circle method to  the second symmetrized rank
moment $\eta_2(n)$,  Bringmann \cite{Bringmann-2008} obtained an asymptotic expression of the spt-function $\spt(n)$.

\begin{thm}[Bringmann]\label{thm-asy-p-o}
As $n\rightarrow\infty$,
\begin{equation}\label{thm-asy-p-o-eqn}
\spt(n)\sim\frac{\sqrt{6}}{\pi}\sqrt{n}p(n)\sim\frac{1}{2\sqrt{2}\pi\sqrt{n}}
e^{\pi\sqrt{\frac{2n}{3}}}.
\end{equation}
\end{thm}

 The above formula also  follows from  an asymptotic estimate of the
  difference of the positive rank moments and the positive  crank moments, due to  Bringmann and Mahlburg \cite{Bringmann-Mahlburg-2014}.

\begin{thm}[Bringmann and Mahlburg]\label{thm-asy-p}
For $r\geq 1$, as $n\rightarrow\infty$,
\[M_r^+(n)-N_r^+(n)\sim \delta_rn^{\frac{r}{2}-\frac{3}{2}}e^{\pi\sqrt{\frac{2n}{3}}},\]
where
\[\delta_r =r!\zeta(r-2)\left(1-2^{3-r}\right)\frac{6^{\frac{r-1}{2}}}{4\sqrt{3}\pi^{r-1}}.\]
\end{thm}

Using Hardy and Ramanujan's   asymptotic formula
\[p(n)\sim\frac{1}{4\sqrt{3}n}e^{\pi\sqrt{\frac{2n}{3}}}, \quad \text{as } \quad n\rightarrow \infty,\]
the $r=2$ case of  Theorem \ref{thm-asy-p} implies Theorem \ref{thm-asy-p-o}, since
\[\spt(n)=M^{+}_2(n)-N^{+}_2(n).\]
Bringmann and Mahlburg \cite{Bringmann-Mahlburg-2014} pointed out that for $r=1$, Theorem \ref{thm-asy-p} leads to an asymptotic formula for   $\ospt(n)$, as  defined in \eqref{defi-ospt-eq}.

 \begin{thm}[Bringmann and Mahlburg]
As $n\rightarrow \infty$,
\begin{equation}\label{BriMaheva-for}
\ospt(n)\sim \frac{p(n)}{4}\sim\frac{1}{16\sqrt{3}n}e^{\pi\sqrt{\frac{2n}{3}}}.
\end{equation}
\end{thm}

 Eichhorn and Hirschhorn \cite{Eichhorn-Hirschhorn-2015} provided  an alternative proof of Theorem \ref{thm-asy-p-o}. In fact, they showed that
 \begin{equation}\label{equ-spt-asy-mean}
 \frac{\spt(n)}{p(n)}\sim\frac{\sqrt{6}}{\pi}\sqrt{n}, \quad \text{as} \quad n\rightarrow \infty.
 \end{equation}
  Let $\lambda$ be a partition of $n$,  define  $n_s(\lambda)$ to be the number of smallest parts of $\lambda$. It is clear that  the left-hand side of \eqref{equ-spt-asy-mean} can be viewed as the mean   of the statistic $n_s(\lambda)$
  over all partitions of $n$.   Eichhorn and Hirschhorn \cite{Eichhorn-Hirschhorn-2015} obtained
   formulas for the  mean and the standard deviation of $n_s(\lambda)$.

\begin{thm}[Eichhorn and Hirschhorn] As $n\rightarrow\infty$, the statistic $n_s(\lambda)$  is distributed roughly as a negative exponential, with mean
\begin{equation}\label{mean-spt}
\mu=\frac{\sqrt{6}}{\pi}\sqrt{n}+\frac{3}{\pi^2}+\text{O}\left(\frac{1}{\sqrt{n}}\right)
\end{equation}
and standard derivation
\begin{equation}
\sigma=\frac{\sqrt{6}}{\pi}\sqrt{n}-\frac{1}{4}+\text{O}\left(\frac{1}{\sqrt{n}}\right).
\end{equation}
\end{thm}

An asymptotic formula with a power saving error term for $\spt(n)$ has been obtained by
Banks, Barquero-Sanchez, Masri, Sheng \cite{Banks-Josiah-Masri-Sheng-2015} based on an asymptotic formula for $p(n)$ due to Masri \cite{Masri-2015}.

In analogy with the explicit formula for $p(n)$ due to Rademacher \cite{Rademacher-1936, Rademacher-1940, Rademacher-1943},  Ahlgren and Andersen \cite{Ahlgren-Scott-Andresen-2016} obtained
an exact expression for the spt-function.

\begin{thm}[Ahlgren and Andersen]
For $n\geq 1$,
\[\spt(n)=\frac{\pi}{6}(24n-1)^{\frac{1}{4}}\sum_{c=1}^\infty \frac{A_c(n)}{c}(I_{1/2}-I_{3/2})\left(\frac{\pi\sqrt{24n-1}}{6c}\right),\]
where $I_\nu$ is the I-Bessel function, $A_c(n)$ is the Kloosterman sum
\[A_c(n)=\sum_{d \mod c \atop (d,c)=1}e^{\pi i s(d,c)-2i\pi\frac{dn}{c}},\]
and $s(d,c)$ is the Dedekind sum
\[s(d,c)=\sum_{r=1}^{c-1}\frac{r}{c}\left(\frac{dr}{c}-\left \lfloor\frac{dr}{c}\right\rfloor-\frac{1}{2}\right).\]
\end{thm}

Asymptotic properties of generalizations and variations of the spt-function have also been well-studied. Recall that the  higher order spt-function $\spt_k(n)$ introduced by Garvan is defined in \eqref{def-sptk}. Its asymptotic property was first conjectured by Bringmann and Mahlburg \cite{Bringmann-Mahlburg-2009}, and then confirmed by Bringmann, Mahlburg and Rhoades   \cite{Bringmann-Mahlburg-Rhoades-2011}.

\begin{thm}[Bringmann, Mahlburg and Rhoades]
As $n\rightarrow\infty$,
\[\spt_k(n)\sim\beta_{2k}n^{k-\frac{1}{2}}p(n),\]
where $\beta_{2k}\in\frac{\sqrt{6}}{\pi}\mathbb{Q}$ is positive.
\end{thm}

The following asymptotic formula for $\Spt_j(n)$, as  defined in \eqref{def-sptj},
  is due to Rhoades \cite{Rhoades-2013}.

\begin{thm}[Rhoades]
As $n\rightarrow \infty$,
\[\Spt_j(n)=\frac{j}{2\pi\sqrt{2n}}e^{\pi\sqrt{\frac{2n}{3}}}
(1+o_j(1)).\]
\end{thm}

 Waldherr \cite{Waldherr-2013} obtained an asymptotic property of  the  $j$-rank moment $_jN_{k}(n)$ defined in \eqref{jrankmom}.

 \begin{thm}[Waldherr]
For $1\leq j\leq 12$,  as $n\rightarrow \infty$,
\begin{equation}
_jN_{2k}(n)\sim 2\sqrt{3}(-1)^{k}B_{2k}\left(\frac{1}{2}\right)(24n)^{k-1}e^{\pi\sqrt{\frac{2n}{3}}},
\end{equation}
where $B_r(\cdot)$ is a Bernoulli polynomial. Furthermore,
\begin{equation}
_{j-1}N_{2k}(n)-{_{j}}N_{2k}(n)\sim \sqrt{3}\frac{(2k)!}{(2k-2)!}(-1)^{k+1}B_{2k-2}(24n)^{k-\frac{3}{2}}
e^{\pi\sqrt{\frac{2n}{3}}}.
\end{equation}
In particular, $_{j-1}N_{2k}(n)>\ _{j}N_{2k}(n)$ for all sufficiently large $n$.
\end{thm}

Kim,  Kim and  Seo  \cite{Kim-Kim-Seo} derived  an asymptotic expression for $\overline{\ospt}(n)$, as defined in \eqref{defi-ospt-b-eqn}.

\begin{thm}[Kim,  Kim and  Seo]\label{Kim-Kim-Seo-thm}
As $n\rightarrow\infty$,
\begin{equation}\label{Kim-Kim-Seo-eqn}
\overline{\ospt}(n)\sim\frac{1}{64n}e^{\pi\sqrt{n}}
\sim \frac{\bar{p}(n)}{8},
\end{equation}
where $\bar{p}(n)$ denotes the number of overpartitions of $n$.
\end{thm}

The above theorem is a  consequence of an asymptotic formula
for the difference of the modified rank and crank moments for overpartitions due to Rolon \cite{Rolon-2016}.

\begin{thm}[Rolon]\label{Rolon-asympt}
As $n\rightarrow\infty$,
\begin{eqnarray}\label{Rolon-asympt-eqn}
\overline{M}_r^+(n)-\overline{N}_r^+(n)
\sim\delta_rn^{\frac{r}{2}-\frac{3}{2}}e^{\pi\sqrt{n}},
\end{eqnarray}
where
\[\delta_r=r!\pi^{-r+1}2^{r-5}\zeta(r-2)\left(1-2^{3-r}\right).\]
\end{thm}

Combining \eqref{defi-ospt-b-eqn} and \eqref{Rolon-asympt-eqn} with $r=1$,
 we arrive at \eqref{Kim-Kim-Seo-eqn}.

\section{Conjectures on inequalities}

In this section, we pose some conjectures on inequalities on the spt-function,
which are reminiscent of inequalities on $p(n)$. We first state some results and conjectures on $p(n)$. Then we
present corresponding conjectures on $\spt(n)$.

 Recall that a sequence $\{a_n\}_{n\geq 0}$ is called log-concave if for  $n\geq 1$,
\begin{equation}\label{log-concave}
a_{n}^2-a_{n-1}a_{n+1}\geq 0.
\end{equation}

It was conjectured in \cite{Chen-2010} that the partition function $p(n)$ is log-concave for $n\geq 26$, that is,
(\ref{log-concave}) is true for  $p(n)$ when $n\geq 26$.  DeSalvo and Pak \cite{DeSalvo-Pak-2015}
confirmed this conjecture by using the Hardy-Ramanujan-Rademacher
formula for $p(n)$ and Lehmer's error bound.

\begin{thm}[DeSalvo and Pak]\label{thm-p-log-c-thm}
For $n\geq 26$,
\begin{equation} \label{thm-p-log-c}
p(n)^2>p(n-1)p(n+1).
\end{equation}
\end{thm}

They also proved the following  inequalities conjectured in  \cite{Chen-2010}.

\begin{thm}[DeSalvo and Pak]
For  $n\geq 2$,
\begin{equation} \label{thm-p-log-chen1}
\frac{p(n-1)}{p(n)}\left(1+\frac{1}{n}\right)>\frac{p(n)}{p(n+1)}.
\end{equation}
\end{thm}

\begin{thm} [DeSalvo and Pak]
For  $n>m>1$,
\begin{equation} \label{thm-p-log-chen2}
p(n)^2\geq p(n-m)p(n+m).
\end{equation}
\end{thm}

DeSalvo and Pak further proved that the term $(1+1/n)$ in \eqref{thm-p-log-chen1}  can be improved  to $(1+O(n^{-3/2}))$.

\begin{thm}[DeSalvo and Pak]
For   $n\geq7$,
\begin{equation} \label{thm-p-log-chen3}
\frac{p(n-1)}{p(n)}\left(1+\frac{240}{(24n)^{3/2}}\right)>\frac{p(n)}{p(n+1)}.
\end{equation}
\end{thm}

DeSalvo and Pak \cite{DeSalvo-Pak-2015} conjectured that the coefficient of $1/n^{3/2}$ in the inequality \eqref{thm-p-log-chen3} can be improved to $\pi/\sqrt{24}$, which was proved by Chen, Wang and Xie \cite{Chen-Wang-Xie-2016}.

\begin{thm}[Chen, Wang and Xie]\label{chen-wang-xie-the}
For   $n\geq 45$,
\begin{equation}\label{Chen-Wang-Xie}
\frac{p(n-1)}{p(n)}\left(1+\frac{\pi}{\sqrt{24}n^{3/2}}\right)>\frac{p(n)}{p(n+1)}.
\end{equation}
\end{thm}

Bessenrodt and Ono \cite{Bessenrodt-Ono-2016} obtained an inequality on $p(n)$.

\begin{thm}[Bessenrodt and Ono]
If $a,b$ are integers with $a,b>1$ and $a+b>8$, then
\begin{equation}\label{BO-ineq}
p(a)p(b)\geq p(a+b),
\end{equation}
where the equality can occur only if  $\{a,b\}=\{2,7\}.$
\end{thm}

We now turn to conjectures on $\spt(n)$.

\begin{conje}\label{conj-logcav-spt-con}
For $n\geq 36$,
\begin{equation} \label{conj-logcav-spt}
\spt(n)^2>\spt(n-1)\spt(n+1).
\end{equation}
\end{conje}

\begin{conje}\label{conj-chen1}
For   $n\geq 13$,
\begin{equation}
\frac{\spt(n-1)}{\spt(n)}\left(1+\frac{1}{n}\right)>\frac{\spt(n)}{\spt(n+1)}.
\end{equation}
\end{conje}

Like the case for $p(n)$, the term $(1+1/n)$ in Conjecture \ref{conj-chen1} can be sharpened to $(1+O(n^{-3/2}))$.

\begin{conje}\label{conj-chen1a}
For   $n\geq 73$,
\begin{equation}\label{spt-inverse}
\frac{\spt(n-1)}{\spt(n)}\left(1+\frac{\pi}{\sqrt{24}n^{3/2}}\right)>\frac{\spt(n)}{\spt(n+1)}.
\end{equation}
\end{conje}

The following conjectures are analogous to  \eqref{thm-p-log-chen2} and \eqref{BO-ineq}.

\begin{conje}
For  $n>m>1$,
\begin{equation} \label{conj-chen2}
\spt(n)^2> \spt(n-m)\spt(n+m).
\end{equation}
\end{conje}

\begin{conje}
If $a,b$ are integers with $a,b>1$  and $(a,b)\neq (2,2)$ or $(3,3)$,  then
\begin{equation}
\spt(a)\spt(b)> \spt(a+b).
\end{equation}
\end{conje}

Beyond quadratic inequalities, we observe that many combinatorial sequences
 including $\{p(n)\}_{n\geq 1}$ and $\{\spt(n)\}_{n\geq 1}$ seem to satisfy higher order inequalities
 except for a few terms at the beginning.
Notice that  $ I(a_0, a_1, a_2)=a_1^2-a_{0}a_{2} $  is an invariant of  the
quadratic binary form
\[a_{2}x^2+2a_1xy+a_{0}y^2.\]
For a sequence $a_0, a_1, a_2, \ldots$ of indeterminates, let
\[ I_{n-1}(a_0, a_1, a_2) = I(a_{n-1}, a_n, a_{n+1}) =a_{n}^2-a_{n-1}a_{n+1}.\]
Then Conjecture \ref{conj-logcav-spt-con} says that for $a_{n}=\spt(n)$, $I_{n-1}(a_0, a_1, a_2)>0$ holds when $n\geq 36$.

This phenomenon occurs for other invariants as well.
For the background on the theory of invariants, see, for example, Hilbert \cite{Hilbert}, Kung and Rota \cite{Kung-Rota-1984} and Sturmfels \cite{Sturmfels}.
A binary form $f(x,y)$ of degree $n$ is a homogeneous polynomial of degree $n$
in two variables $x$ and $y$:
 \[f(x,y)=\sum_{i=0}^n {n\choose i} a_ix^iy^{n-i},\]
where the coefficients $a_i$ are complex numbers.

Let
\[C=\left(\begin{array}{ll}
c_{11} &c_{12}\\
c_{21} &c_{22}\end{array}\right)\]
be an invertible complex matrix. Under the  linear transformation
\begin{align*}
x&=c_{11} \overline{x}+c_{12} \overline{y},\\[5pt]
y&=c_{21}\overline{x}+c_{22} \overline{y},
\end{align*}
 the binary form  $f(x,y)$ is transformed into another binary form
\begin{align*}
\overline{f} (\overline{x},\overline{y})
&=
\sum_{i=0}^n {n\choose i} \overline{a}_{i}\, \overline{x}^i\, \overline{y}^{n-i},
\end{align*}
where the coefficients $ \overline{a}_i$ are   polynomials in  $a_i$ and $c_{ij}$.
Let  $g$ be a nonnegative integer. A polynomial  $I({a}_{0},{a}_{1},\ldots,{a}_{n})$ in the coefficients $a_0,a_1,\ldots,a_n$ is    an invariant of index $g$ of the binary form $f(x,y)$ if  for any invertible matrix $C$,
\[
  I({\overline{a}}_{0},{\overline{a}}_{1},\ldots,{\overline{a}}_{n})
  = {({c}_{11}{c}_{22}-{c}_{12}{c}_{21})}^{g}I({a}_{0},{a}_{1},\ldots,{a}_{n}).
\]
For example,
\begin{equation}\label{equ-dis-cub}
  I(a_0, a_1, a_2, a_3) = 3a_1^2 a_2^2-4a_1^3 a_3-4a_0 a_2^3-a_0^2 a_3^2+6a_0a_1a_2a_3
  \end{equation}
is an invariant of the  cubic binary  form
\begin{equation}\label{cubicform}
f(x,y)=a_3x^3+3a_2x^2y+3a_1xy^2+a_0y^3.
\end{equation}
Note that $27I(a_0, a_1, a_2, a_3)$ is called the discriminant of \eqref{cubicform}.
The polynomial $I(a_{n-1}, a_n, a_{n+1}, a_{n+2})$ is related to  the higher order Tur\'an inequality.
Recall that a sequence $\{a_n\}_{ n\geq 0}$ satisfies the higher order Tur\'{a}n inequality   if for $ n\geq 1$,
\begin{equation}\label{high-turan}
4(a_n^2-a_{n-1}a_{n+1})(a_{n+1}^2-a_{n}a_{n+2})-(a_na_{n+1}-a_{n-1}a_{n+2})^2> 0,
\end{equation}
and we say that  $\{a_n\}_{ n\geq 0}$ satisfies the Tur\'{a}n inequality if it is log-concave.

A simple calculation
shows that for
 $n=1$, the polynomial in \eqref{high-turan}
reduces to the invariant $I(a_0, a_1, a_2, a_3)$ in \eqref{equ-dis-cub}, namely,
\begin{align*}
&\!\!\!\!  3a_1^2 a_2^2-4a_1^3 a_3-4a_0 a_2^3-a_0^2 a_3^2+6a_0a_1a_2a_3\\[3pt]
&\qquad= 4(a_1^2-a_{0}a_{2})(a_{2}^2-a_{1}a_{3})-(a_1a_{2}-a_{0}a_{3})^2.
\end{align*}

Csordas, Norfolk and Varga \cite{Csordas-Norfolk-Varga-1986} proved
that  the coefficients of the Riemann $\xi$-function
satisfy the Tur\'{a}n inequality. This settles a conjecture of  P\'{o}lya.
Dimitrov \cite{Dimitrov-1998} showed   under the
Riemann hypothesis, the coefficients of the Riemann $\xi$-function
satisfy the higher order Tur\'{a}n inequality. Dimitrov and Lucas \cite{Dimitrov-Lucas-2011}
proved this assertion without  the Riemann hypothesis.

Numerical evidence indicates that
both  $p(n)$ and  $\spt(n)$ satisfy  the high order Tur\'an inequality.

 \begin{conje}\label{conj-dis-3}
For $n\geq 95$,  $p(n)$ satisfies the higher order Tur\'an
 inequality \eqref{high-turan},  whereas   $\spt(n)$
satisfies  \eqref{high-turan} for $n\geq 108$.
\end{conje}

We next consider the invariant of the
quartic binary form
\begin{equation}\label{equ-quartic-bi-f}
f(x,y)=a_4x^4+4a_3x^3y+6a_2x^2y^2+4a_1xy^3+a_0y^4.
\end{equation}
It appears that for large $n$, both $p(n)$ and  $\spt(n)$ satisfy
the inequalities derived from  the following invariants of  \eqref{equ-quartic-bi-f}:
  \begin{eqnarray*}
   A(a_0,a_1,a_2,a_3,a_4)&=&a_0a_{4}-4a_{1}a_{3}+3a_{2}^2 ,\\[5pt]
   B(a_0,a_1,a_2,a_3,a_4)&=&-a_0a_{2}a_{4}+a_{2}^3+a_0a_{3}^2
    +a_{1}^2a_{4}-2a_{1}a_{2}a_{3},\\[5pt]
    I(a_0,a_1,a_2,a_3,a_4)&=&A(a_0,a_1,a_2,a_3,a_4)^3-27B(a_0,a_1,a_2,a_3,a_4)^2.
  \end{eqnarray*}
Notice that $256I(a_0,a_1,a_2,a_3,a_4)$ is the discriminant of $f(x,y)$ in \eqref{equ-quartic-bi-f}.
  To be more specific, we have the following conjectures:
Setting $a_n=p(n)$,
\begin{eqnarray}
A(a_{n-1},a_n,a_{n+1},a_{n+2},a_{n+3})&>&0,\qquad\text{for }n\geq 185,\label{conj-partition-A}\\[3pt]
B(a_{n-1},a_n,a_{n+1},a_{n+2},a_{n+3})&>&0,\qquad\text{for }n\geq 221,\label{conj-partition-B}\\[3pt]
I(a_{n-1},a_n,a_{n+1},a_{n+2},a_{n+3})&>&0,\qquad\text{for }n\geq 207.
\end{eqnarray}
Setting  $a_n=\spt(n)$,
\begin{eqnarray}
A(a_{n-1},a_n,a_{n+1},a_{n+2},a_{n+3})&>&0,\qquad\text{for }n\geq 205,\label{conj-spt-A}\\[3pt]
B(a_{n-1},a_n,a_{n+1},a_{n+2},a_{n+3})&>&0,\qquad\text{for }n\geq 241,
\label{conj-spt-B}\\[3pt]
I(a_{n-1},a_n,a_{n+1},a_{n+2},a_{n+3})&>&0,\qquad\text{for }n\geq 227.
\end{eqnarray}

In general, it would be interesting to further study higher order
inequalities on $p(n)$ and $\spt(n)$ based on polynomials arising in
the invariant theory of binary forms.

After the submission of an early version for the proceedings of the
26th British Combinatorial Conference (Surveys in Combinatorics 2017,
 A. Claesson, M. Dukes, S. Kitaev, D. Manlove and K. Meeks, Eds.,  Cambridge
 University Press, Cambridge, 2017), we observed that the above conjectured
inequalities on $p(n)$ and $\spt(n)$ seem to permit companion inequalities analogous to \eqref{Chen-Wang-Xie} and \eqref{spt-inverse}.

\begin{conje}\label{conj-chen-inverse-higher-partition}
Let
\begin{equation}\label{partition-u}
u_n=\frac{p(n+1)p(n-1)}{p(n)^2}.
\end{equation}
For $n\geq 2$,
\begin{equation}
\left(1+\frac{\pi}{\sqrt{24}n^{3/2}}\right)\left(1-u_n u_{n+1}\right)^2>4\left(1-u_n\right)\left(1-u_{n+1}\right).
\end{equation}
\end{conje}

\begin{conje}\label{conj-chen-inverse-higher-spt}
Let
\begin{equation}\label{spt-v}
v_n=\frac{\spt(n+1)\spt(n-1)}{\spt(n)^2}.
\end{equation}
For $n\geq 2$,
    \begin{equation}
\left(1+\frac{\pi}{\sqrt{24}n^{3/2}}\right)\left(1-v_n v_{n+1}\right)^2
>4\left(1-v_n\right)\left(1-v_{n+1}\right).
    \end{equation}
\end{conje}

As for the inequalities
\eqref{conj-partition-A}, \eqref{conj-partition-B}, \eqref{conj-spt-A} and \eqref{conj-spt-B} for $p(n)$ and $\spt(n)$ based on the
invariants $A(a_0,a_1,a_2,a_3,a_4)$ and $B(a_0,a_1,a_2,a_3,a_4)$
of the quartic binary form \eqref{equ-quartic-bi-f},
it appears that there exist similar companion inequalities.

\begin{conje}
We have
\begin{equation}
4\Big(1+\frac{\pi^2}{16n^3}\Big){ a_{n}a_{n+2}}>{a_{n-1}a_{n+3}+3a_{n+1}^2}
\end{equation}
for   $a_n=p(n)$ when $n>217$ and  for $a_n=\spt(n)$ when
 $n>259$.
\end{conje}

\begin{conje}
We have
\begin{align}
\nonumber &\Big(1+\frac{\pi^3}{72\sqrt{6}n^{9/2}}\Big)\big(
{2a_{n}a_{n+1}a_{n+2}+a_{n-1}a_{n+1}a_{n+3}}\big)&\\[3pt]
 &\qquad \qquad > {a_{n+1}^3+a_{n-1}a_{n+2}^2+a_{n}^2a_{n+3}}
\end{align}
 for $a_n=p(n)$ when  $n>243$ and for $a_n=\spt(n)$ when $n>289$.
\end{conje}

\vskip 0.2cm

 \noindent{\bf Acknowledgments.} We wish to thank  George   Andrews, Kathrin  Bringmann, Mike Hirschhorn, Joseph Kung,  Karl Mahlburg, Ken Ono and Peter Paule for valuable comments and suggestions. This work was supported by the National Science Foundation of China.

%\newpage
%\thispagestyle{empty}
\mbox{}
%\newpage

\end{document}